\documentclass{article}
\usepackage{amsmath,amssymb,amsfonts,amscd}
\usepackage{hyperref}
\usepackage[numeric]{amsrefs}

\newtheorem{observation}[equation]{Observation}
\newtheorem{thm}[equation]{Theorem}
\newtheorem{prop}[equation]{Proposition}
\newtheorem{lem}[equation]{Lemma}

\DeclareMathOperator{\sg}{sgn}

\numberwithin{equation}{section}
\renewcommand\a{\alpha}
\renewcommand\b{\beta}
\newcommand\g{\gamma}
\renewcommand\d{\delta}
\newcommand\e{\varepsilon}
\renewcommand\l{\lambda}
\renewcommand\L{\Lambda}
\newcommand\G{\Gamma}
\newcommand\f{\frac}

\newcommand{\N}{{\mathbb{N}}}
\newcommand{\Z}{{\mathbb{Z}}}
\newcommand{\R}{{\mathbb{R}}}
\newcommand{\RP}{{\mathbb{RP}}}
\newcommand{\C}{{\mathbb{C}}}

\newcommand{\Q}{{\mathbb{Q}}}

\newcommand\re{\text{Re~}}

\renewcommand\Re{\text{Re~}}
\renewcommand\Im{\mbox{Im~}}

\newcommand{\tr}{\operatorname{tr}}

\newcommand{\sgn}{\operatorname{sgn}}

\renewcommand\H{{\cal H}}

\renewcommand\i{^{-1}}
\renewcommand\({\left(}
\renewcommand\){\right)}
\newcommand{\ttwo}[4]{
\(\begin{smallmatrix}{#1} & {#2}
\\ {#3} & {#4} \end{smallmatrix}\)}
\newcommand{\bx}{\hfill$\square$\vspace{.6cm}}

\newcommand{\gobble}[1]{}
  \newcommand{\rangeref}[2]{%
    \ref{#1}--\afterassignment\gobble\fam 0\ref{#2}%
  }

\begin{document}

\title{The Rankin-Selberg method for automorphic distributions}
\author{Stephen D. Miller\thanks{Partially supported
by NSF grant DMS-0301172 and an Alfred P. Sloan Foundation
Fellowship.}~ and  Wilfried Schmid\thanks{Partially supported by
NSF grant DMS-0500922.}}
\date{December 7, 2005}

\maketitle

\renewcommand{\theequation}{\thesection.\arabic{equation}}%

\section{Introduction}\label{intro}

We recently established the holomorphic continuation and functional equation of the exterior square $L$-function for
$GL(n,\Z)$, and more generally, the archimedean theory of the $GL(n)$ exterior square $L$-function over $\Q$. We
refer the reader to our paper \cite{extsquare} for a precise statement of the results and their relation to previous work
on the subject. The purpose of this note is to give an account of our method in the simplest non-trivial cases, which
can be explained without the technical overhead necessary for the general case.

Let us begin by recalling the classical results, about standard $L$-functions and Rankin-Selberg $L$-functions
of modular forms. We consider a cuspidal modular form $F$, of weight $k$, on the upper half plane $\H$. To simplify
the notation, we suppose that it is automorphic for $\G=SL(2,\Z)$, though the arguments can be adapted to congruence subgroups of $SL(2,\Z)$. Like all modular forms, $F$ has a Fourier expansion,
\begin{equation}
\label{modular1}
F(z)\ = \ {\sum}_{n\geq 1}\, a_n\,e(nz)\,,\ \ \ \text{with}\ \ e(z)\ =_{\text{def}}\ e^{2\pi i z}\,.
\end{equation}
For a general modular form, the Fourier series may involve a non-zero constant term $a_0$; it is the hypothesis
of cuspidality that excludes the constant term. The Dirichlet series
\begin{equation}
\label{modular2}
L(s,F)\ = \ {\sum}_{n\geq 1}\, a_n\,n^{-\frac{k-1}{2}-s}
\end{equation}
converges for $\Re s \gg 0$, extends holomorphically to the entire $s$-plane, and satisfies a functional equation
relating $L(s,F)$ to $L(1-s,F)$. This is the standard $L$-function of the modular form $F$.

Hecke proved the holomorphic continuation and functional equation by expressing $L(s,F)$ in terms
of the Mellin transform of $F$ along the imaginary axis,
\begin{equation}
\label{modular3}
\begin{aligned}
&\int_0^\infty F(iy)\,y^{s-1}\,dy\ = \ {\sum}_{n\geq 1}\, a_n \int_0^\infty e^{-2\pi n y} \,y^{s-1}\,dy
\\
&\ \ \ \,  = \, {\sum}_{n\geq 1}\, a_n \,n^{-s} \int_0^\infty\! e^{-2\pi y} \,y^{s-1}\,dy \, = \, (2\pi)^{-s}\,\G(s)\, L(s - \textstyle\frac{k-1}{2},F)\,,
\end{aligned}
\end{equation}
at least for $\Re s \gg 0$. The transformation law for the modular form $F$ under $z\mapsto -1/z$,
\begin{equation}
\label{modular4}
F(-1/z)\ = \ (-z)^{k}\,F(z)\,,
\end{equation}
implies that $F(iy)$ decays rapidly not only as $y\to\infty$, but also as $y\to 0$. That makes the Mellin transform, and
hence also $\G(s+\frac{k-1}2)\, L(s,F)$, globally defined and holomorphic. The Gamma function has no zeroes, so
$L(s,F)$ is entire as well. The transformation law (\ref{modular4}), coupled with the change of variables $y \mapsto 1/y$
and the shift $s\mapsto s+\frac{k-1}2$, gives the functional
equation
\begin{equation}
\label{modular5}
\begin{aligned}
&(2\pi)^{-s-\frac{k-1}2}\,\G(s+\textstyle\frac{k-1}2)\, L(s,F)\ =
\\
&\qquad\qquad\qquad = \  i^k\,(2\pi)^{s-1-\frac{k-1}2}\,\G(1-s+\textstyle\frac{k-1}2)\, L(1-s,F)\,.
\end{aligned}
\end{equation}
The factor $i^k$ comes up naturally in the computation, yet might be misleading since $\G=SL(2,\Z)$ admits only
modular forms of even weights.

In addition to $F$, we now consider a second modular form of weight $k$, which need not be cuspidal,
\begin{equation}
\label{modular6}
G(z)\ = \ {\sum}_{n\geq 0}\, b_n\,e(nz)\,.
\end{equation}
The Rankin-Selberg $L$-function of the pair $F$, $\overline G=$ complex conjugate of $G$, is the Dirichlet series
\begin{equation}
\label{modular7}
L(s,F\otimes\overline G)\ = \ \zeta(2s)\,\,{\sum}_{n\geq 1}\, a_n\,\overline b_n\,n^{1-k-s}\,.
\end{equation}
Its analytic continuation and functional equation were established separately by Rankin \cite{rankin} and Selberg
\cite{selberg}. The proof depends on properties of the non-holomorphic Eisenstein series
\begin{equation}
\label{modular8}
\begin{aligned}
E_s(z)\ = \ \ &\pi^{-s}\,\G(s)\,\zeta(2s)\,\,{\sum}_{\gamma\in\G_\infty\backslash\G}\, \bigl(\Im(\gamma z)\bigr)^s
\\
&\bigl(\, \G_\infty \ = \ \{\,\gamma\in\G\,\mid\,\gamma\infty=\infty\,\}\,\bigr).
\end{aligned}
\end{equation}
This sum is well defined since $\G_\infty$ acts on $\H$ by integral translations.  It converges for $\Re s > 1$ and
extends meromorphically to the entire $s$-plane with only one pole, of first order, at $s=1$. The function $E_s(z)$ is
$\G$-invariant by construction, has moderate growth as $\Im z\to\infty$, and satisfies the functional equation
\begin{equation}
\label{modular9}
E_s(z)\ = \ E_{1-s}(z)\,.
\end{equation}
Both $F(z)$ and $G(z)$ transform according to a factor of automorphy under the action of $\G$, but
$(\Im z)^k F(z)\overline G(z)$ is $\G$-invariant, as is the measure $y^{-2}dxdy$. Since $G(z)$ and $E_s(z)$ have
moderate growth as $\Im z\to\infty$, and since $F(z)$ decays rapidly, the integral
\begin{equation}
\label{modular10}
I(s)\ = \ \int_{\G\backslash\H}(\Im z)^{k-2} F(z)\,\overline G(z)\,E_s(z)\,dx\,dy
\end{equation}
converges. From $E_s(z)$, the function $I(s)$ inherits both the functional equation
\begin{equation}
\label{modular11}
I(s)\ = \ I(1-s)
\end{equation}
and the analytic properties: it is holomorphic, with the exception of a potential first order pole at $s=1$.

The definition (\ref{modular8}) of $E_s(z)$ involves a sum of $\,\Im \gamma z$, with $\gamma$ ranging over $\G_\infty\backslash\G$.
But the rest of the integrand in (\ref{modular10}) is $\Gamma$-invariant. That justifies the process known as ``unfolding",
\begin{equation}
\label{modular12}
\begin{aligned}
&\pi^{s}\bigl(\G(s)\,\zeta(2s)\bigr)^{-1}I(s)\ =
\\
&\qquad= \ \int_{\G\backslash \H} {\sum}_{\gamma\in\G_\infty\backslash\G}\,(\Im z)^{k-2} F(z)\,\overline G(z)\,\bigl(\Im(\gamma z)\bigr)^{s} dx\,dy
\\
&\qquad= \ \int_{\G_\infty\backslash \H} \,(\Im z)^{s+k-2} F(z)\,\overline G(z)\,dx\,dy \,,
\end{aligned}
\end{equation}
at least for $\Re s > 1$, in which case the integral on the right converges. Since $\G_\infty$ acts on $\H$ by integral
translations, the strip $\,\{\,0 \leq \Re z \leq 1\,\}\,$ constitutes a fundamental domain for this action. Substituting the series
(\ref{modular1},\,\ref{modular6}) for $F(z)$ and $G(z)$, one finds
\begin{equation}
\label{modular13}
\begin{aligned}
&\pi^{s}\bigl(\G(s)\,\zeta(2s)\bigr)^{-1}I(s)\ =
\\
&\qquad= \ \int_{0}^{\infty}  \! \int_{0}^{1} {\sum}_{\stackrel{\scriptstyle{\!\! n> 0}}{m\geq 0}} \,\,a_n\,\overline b_m\,e\bigl((n-m)x\bigr)\,e^{-2\pi(n+m)y}y^{s+k-2} dx\,dy
\\
&\qquad= \ \ {\sum}_{n\geq 1} \,\, a_n\,\overline b_n \, \int_{0}^{\infty} \,e^{-4\pi ny}\,y^{s+k-2} \,dy
\\
&\qquad= \,\ (4\pi)^{-s-k+1}\,\G(s+k-1)\,\,{\sum}_{n\geq 1} \, a_n\,\overline b_n \, n^{-s-k+1} \,,
\end{aligned}
\end{equation}
again for $\Re s > 1$. Equivalently,
\begin{equation}
\label{modular14}
I(s)\ = \ 2^{1-k}\,(2\pi)^{1-k-2s}\,\G(s)\,\G(s+k-1)\, L(s,F\otimes \overline G)\,.
\end{equation}
The Gamma factors have no zeroes, so $L(s,F\otimes \overline G)$ extends holomorphically to all of $\C$, except possibly
for a first order pole at $s=1$. In effect, (\ref{modular11}) is the functional equation for the Rankin-Selberg $L$-function. With
some additional effort one can modify these arguments, to make them work even when $F$ and $G$ have
different weights.

Maass \cite{maass} extended the proofs of the analytic continuation and functional equation for the standard $L$-function
to the case of Maass forms, i.e., $\G$-invariant eigenfunctions of the hyperbolic Laplacian on $\H$; see section 2 below.
Jacquet \cite{jacquetsequel} treats the Rankin-Selberg $L$-function for Maass forms. We just saw how the Gamma factors in
(\ref{modular3}) and (\ref{modular13}) arise directly from the standard integral representation of the Gamma function. In
contrast, for Maass forms, the Gamma factor for the standard $L$-function arises from the Mellin transform of the Bessel
function $K_\nu (y)$,
\begin{equation}
\label{modular15}
\int_0^\infty K_\nu(y)\,y^{s-1}\ = \ 2^{s-2}\,\G({\textstyle\frac{s-\nu}{2}})\,\G({\textstyle\frac{s+\nu}{2}})\qquad (\,\Re s \gg 0\,)\,,
\end{equation}
and for the Rankin-Selberg $L$-function of a pair of Maass forms, from the integral
\begin{equation}
\label{modular16}
\begin{aligned}
&\!\!\!\int_0^\infty K_\mu(y)\,K_\nu(y)\,y^{s-1}\,dy\ =
\\
&\ \ \  = \ 2^{s-3}\,\f{\G(\frac{s -  \mu  - \nu }{2})\, \G(\frac{s + \mu  - \nu }{2})\,  \G(\frac{s - \mu  + \nu }{2})\,
\G(\frac{s + \mu  + \nu }{2})}{\G(s)} \ \ \ \ (\, \Re{s} \gg 0 \,)\,.
\end{aligned}
\end{equation}
Though (\ref{modular16}) can be established by elementary means, it is still complicated and its proof lacks a
conceptual explanation.

In the case of Rankin-Selberg $L$-functions of higher rank groups, the integrals analogous to (\ref{modular16})
become exceedingly difficult, or even impossible, to compute. In fact, it is commonly believed that such integrals
may not always be expressible in terms of Gamma functions \cite{bump}*{\S 2.6}. If true, this would not contradict
Langlands' prediction that the functional equations involve certain definite Gamma factors \cites{langlandsdc,
eulerproducts}: the functional equations pin down only the ratios of the Gamma factors on the two sides, which
can of course be expressed also as ratios of other functions.

Broadly speaking, the existing approaches to the $L$-functions for higher rank groups overcome the problem
of computing these so-called {\it archi\-me\-dean integrals} in one of two ways. Even if the integrals cannot be
computed explicitly, it may be possible to establish a functional equation with unknown coefficients; it may then
be possible to identify the coefficients in some special case, or by an analysis of their zeroes and poles. The
Langlands-Shahidi method, on the other hand, often exhibits the functional equation with precisely the Gamma
factors predicted by Langlands. Both methods have one difficulty in common: ruling out poles -- other than those
at the expected places -- of the $L$-functions in question requires considerable effort, and is not always possible.

We are approaching the analytic continuation and functional equation of $L$-functions from a different point of view.
Instead of working with automorphic forms -- i.e., the higher dimensional analogues of modular forms and Maass
forms -- we attach the $L$-functions to {\it automorphic distributions}. In the case of modular forms and Maass forms,
the automorphic distributions can be described quite concretely as boundary values. Alternatively but equivalently,
they can be described abstractly; see \cite{voronoi}*{\S 2} or section three below. Computing with distributions presents
some technical difficulties. What we gain in return are explicit formulas for the archimedean integrals that arise in the
setting of automorphic distributions. This has led us to some new results.

In the next section we show how our method works in the simplest case, for the standard $L$-functions of modular
forms and Maass forms. We treat the Rankin-Selberg $L$-function in section four, following the description of our
main analytic tool in section three. Section five, finally, is devoted to the exterior square $L$-function for $GL(4,\Z)$.
That is the first not-entirely-trivial case of the main result of \cite{extsquare}. It can be explained more transparently
than the general case for two reasons: the main analytic tool is the pairing of distributions, which for $GL(4)$ reduces
to a variant of the Rankin-Selberg method for $GL(2)$. Also, the general case involves a somewhat subtle induction,
with $GL(4)$ representing merely the initial step.

\section{Standard $L$-functions for $SL(2)$}\label{GL2}

Holomorphic functions on the disk or the upper half plane have hyperfunction boundary values, essentially by
definition of the notion of hyperfunction. Holomorphic functions of moderate growth, in particular modular
forms, have distribution boundary values:
\begin{equation}
\label{gl2,1}
\tau(x)\ = \ {\lim}_{y\to 0^+}\, F(x+iy)
\end{equation}
is the automorphic distribution corresponding to a modular form $F$ for $SL(2,\Z)$, of weight $k$. The limit
exists in the strong distribution topology. From $F$, the distribution $\tau$ inherits its $SL(2,\Z)$-automorphy property
\begin{equation}
\label{gl2,2}
\tau(x)\ = \ (cx + d)^{-k}\,\tau\bigl({\textstyle\frac{ax+b}{cx+d}}\bigr)\ \ \ \ \text{for all}\ \ \ttwo{a}{b}{c}{d} \in SL(2,\Z)\,.
\end{equation}
In terms of Fourier expansion (\ref{modular1}) of the cuspidal modular form $F(z)$, the limit (\ref{gl2,1}) can be taken term-by-term,
\begin{equation}
\label{gl2,4}
\tau(x)\ = \ {\sum}_{n>0}\, a_n\,e(nx)\,.
\end{equation}
We shall argue next that it makes sense to take the Mellin transform of the distribution $\tau$, and that this Mellin
is an entire function of the variable $s$. The argument will be a special case of the techniques developed in our
paper \cite{inforder}.

Note that the periodic distribution $\tau$ has no constant term. It can therefore be expressed as the $\ell$-th derivative of a
continuous, periodic function $\phi_\ell$, for every sufficiently large integer $\ell$,
\begin{equation}
\begin{aligned}
\label{gl2,5}
&\tau(x)\ = \ \phi_\ell ^{(\ell)}(x)\,, \ \ \text{with}\ \ \phi_\ell \in  C(\R/\Z)
\\
&\qquad\qquad\bigl(\,\phi_\ell (x)\ = \ {\sum}_{n>0}\, (2\pi i n)^{-\ell}\,a_n\,e(nx)\,\bigr).
\end{aligned}
\end{equation}
Using the formal rule for pairing the ``test function" $\,x^{s-1}\,$ against the deri\-vative of a distribution, we find
\begin{equation}
\label{gl2,6}
\int_0^\infty\!\!\!\! x^{s-1}\,\tau(x)\,dx\, = \int_0^\infty\!\!\!\! x^{s-1}{\textstyle\frac{d^\ell\ }{dx^\ell}}\,\phi_\ell (x)\,dx\, = \, (-1)^\ell \!\! \int_0^\infty\!\!\!\! \phi_\ell (x)\,{\textstyle\frac{d^\ell\ }{dx^\ell}}\,x^{s-1}\,dx \,.
\end{equation}
As a continuous, periodic function, $\phi_\ell$ is bounded. That makes the expression on the right in (\ref{gl2,6})
integrable away from $x=0$, provided $\ell > \Re s$. Indeed, if we multiply the Mellin kernel $x^{s-1}$ by a
cutoff function $\psi\in C^\infty(\R)$, with $\psi(x)\equiv 1$ near $x=\infty$ and $\psi(x)\equiv 0$ near $x=0$, the
resulting integral is an entire function of the variable $s$ -- we simply choose $\ell$ larger than the real part of any
particular $s$. Increasing the value of $\ell$ further does not affect the integral, as can be seen by a legitimate
application of integration by parts. The identity (\ref{gl2,2}), with $a=d=0$, $b=-c=1$, gives
\begin{equation}
\label{gl2,7}
\tau(x)\ = \ (-x)^{-k}\,\tau(-1/x)\,,
\end{equation}
so the behavior of $\tau(x)$ near zero duplicates its behavior near $\infty$, except for the factor $(-x)^k$ which
can be absorbed into the Mellin kernel. The expression on the right in (\ref{gl2,6}) is therefore integrable even down to
zero, and
\begin{equation}
\label{gl2,8}
s\ \mapsto \int_0^\infty \! \tau(x)\,x^{s-1}\,dx\ \ \text{is a well defined, entire holomorphic function.}
\end{equation}
The change of variables $x\mapsto 1/x$ and the transformation law (\ref{gl2,7}) imply
\begin{equation}
\label{gl2,9}
\int_0^\infty \! \tau(x)\,x^{s-1}\,dx\ =\ (-1)^k\int_0^\infty \! \tau(-x)\,x^{k-s-1}\,dx\,.
\end{equation}
The integral on the right is of course well defined, for the same reason as the integral (\ref{gl2,8}).

In view of the argument we just sketched, it is entirely legitimate to replace $\tau(x)$ by its Fourier series and to
interchange the order of integration and summation: for $\Re s \gg 0$,
\begin{equation}
\begin{aligned}
\label{gl2,10}
&\!\!\!\int_0^\infty \! \tau(x)\,x^{s-1}\,dx\ = \  \int_0^\infty \, {\sum}_{n>0}\,a_n\,e(nx)\,x^{s-1}\,dx
\\
&\ \  = \  {\sum}_{n>0}\,\,  a_n \! \int_0^\infty\!\! e(nx)\,x^{s-1}\,dx\ = \ L(s - {\textstyle\frac{k-1}{2}},F) \! \int_0^\infty \!\! e(x)\,x^{s-1}\,dx\,;
\end{aligned}
\end{equation}
recall (\ref{modular2}). The integral on the right makes sense for $\Re s > 0$ if one regards $e(x)$ as a distribution
and applies integration by parts, as was done in the case of $\tau(x)$. In the range $0 < \Re s < 1$ it converges
conditionally. This integral is well known,
\begin{equation}
\label{gl2,11}
\int_0^\infty \! e(x)\,x^{s-1}\,dx\ =\ (2\pi)^{-s}\,\G(s)\,e(s/4)\qquad (\,0 < \Re s < 1\,).
\end{equation}
Since $\G(s)e(s/4)$ has no zeroes, (\ref{gl2,8}) and (\rangeref{gl2,10}{gl2,11}) imply that $L(s,F)$ is entire. Replacing $\tau(x)$ by
$\tau(-x)$ in (\ref{gl2,10}) has the effect of replacing $e(x)$ by $e(-x)$, and accordingly the factor $e(s/4)$ by
$e(-s/4)$ in (\ref{gl2,11}). Thus (\rangeref{gl2,8}{gl2,11}) imply
\begin{equation}
\label{gl2,12}
\begin{aligned}
&(2\pi)^{-s}\,e(s/4)\,\G(s)\,L(s-{\textstyle\frac{k-1}{2}},F)\ =
\\
&\qquad\qquad =\ (-1)^k\,(2\pi)^{s-k}\, e\bigl((s-k)/4\bigr)\,\G(k-s)\,L(1-s+{\textstyle\frac{k+1}{2}},F)\,.
\end{aligned}
\end{equation}
Since $e(-k/4)=i^{-k}$, this functional equation is equivalent to the functional equation stated in (\ref{modular5}).

A Maass form is a $\G$-invariant eigenfunction $F\in C^\infty(\H)$ for the hyperbolic Laplacian $\Delta$, of moderate growth
towards the boundary of $\H$. It is convenient to express the eigenvalue as $\,(\l^2-1)/4$, so that
\begin{equation}
\label{maass1}
y^2\left({\textstyle\frac{\partial^2\ }{\partial x^2} + \frac{\partial^2\ }{\partial y^2}}\right)F\ = \ {\textstyle\frac {\l^2-1}4}\,F\,.
\end{equation}
Near the real axis, the Maass form $F$ has an asymptotic expansion,
\begin{equation}
\label{maass2}
F(x+iy)\ \sim \ y^{\frac{1-\l}2}\,{\sum}_{k\geq 0}\,\tau_{\l,k}(x)\,y^{2k} \ + \ y^{\frac{1+\l}2}\,{\sum}_{k\geq 0}\,\tau_{-\l,k}(x)\,y^{2k}
\end{equation}
as $y$ tends to zero from above, with distribution coefficients $\,\tau_{\pm\l,k}\,$. In the exceptional case $\,\l=0$, the
leading terms $\,y^{(1-\l)/2}$, $\,y^{(1+\l)/2}\,$ must be replaced by, respectively, $y^{1/2}$ and $y^{1/2}\log y$. The leading
coefficients
\begin{equation}
\label{maass3}
\tau_\l\ =_{\text{def}}\ \tau_{\l,0}\,,\ \ \ \ \tau_{-\l}\ =_{\text{def}}\ \tau_{-\l,0}
\end{equation}
determine the others recursively. They are the {\it automorphic distributions} corresponding to the Maass form $F$. Each
of the two also determines the other -- in a way we shall explain later -- unless $\l$ is a negative odd integer, in which
case the $\tau_{-\l,k}$ all vanish identically. To avoid making statements with trivial counterexamples, we shall not
consider $\,\tau_{-\l}$ when $\l\in \Z_{<0}\cap (2\Z+1)$, and for $\l=0$, we shall only consider the coefficient of
$y^{1/2}$, not the coefficient of $y^{1/2}\log y$.

Unlike modular forms, Maass forms are $\,\G$-invariant as functions, i.e., without a factor of automorphy. However, because
of the nature of the asymptotic expansion (\ref{maass2}), the $\,\G$-invariance of $F$ translates into an automorphy
condition on the automorphic distributions,
\begin{equation}
\label{maass4}
\tau_\l(x)\ =\ |cx+d|^{\l-1}\,\tau_{\l}\left({\textstyle\frac{ax+b}{cx+d}}\right)\ \ \ \text{for all}\ \ \ttwo{a}{b}{c}{d} \in \G\,.
\end{equation}
To simplify the discussion, we suppose $\,\G=SL(2,\Z)$, as before. Then (\ref{maass4}), with $a=b=d=1$, $c=0$ implies
$\tau_\l(x)\equiv \tau_\l(x+1)$, so $\,\tau_\l$ has a Fourier expansion
\begin{equation}
\label{maass5}
\tau_\l(x)\ =\ {\sum}_{n\in\Z}\,a_n\,e(nx)\,.
\end{equation}
From the point of view of $L$-functions, cuspidal Maass forms are more interesting than non-cuspidal forms. The
condition of cuspidality on $F$ is equivalent to two conditions on the automorphic distribution $\,\tau_\l$, namely
\begin{equation}
\label{maass6}
a_0\ = \ 0\,,\ \ \text{and $\,\tau_\l\,$ vanishes to infinite order at $x=0$\,}
\end{equation}
\cite{inforder}. To explain the meaning of the second condition, we note that the discussion leading up to (\ref{gl2,5})
applies also in the present context, since $a_0=0$. The automorphy condition (\ref{maass4}), with $a=d=0$, $b=-c=-1$,
asserts
\begin{equation}
\label{maass7}
\tau_\l(x)\ = \ |x|^{\l -1}\,\tau_\l(-1/x)\,.
\end{equation}
Combined with (\ref{gl2,5}) and the chain rule for the change of variables $x\mapsto -1/x$, this implies
\begin{equation}
\label{maass8}
\tau_\l(x)\ = \ |x|^{\l -1}\,\left(x^2\,{\textstyle\frac{d\ }{dx}}\right)^\ell \bigl(\phi_{\ell}(-1/x)\bigr)\ \ \text{on}\,\ \R-\{0\}\,,
\end{equation}
for every sufficiently large $\ell\in\N$, with some $\phi_\ell\in C(\R-\{0\})$ which remains bounded as $|x|\to \infty$. Moving
the factor $|x|^{\l-1}$ across the differential operator and keeping track of the powers of $x$ shows that the right hand
side of (\ref{maass8}) defines a distribution even on a neighborhood of the origin -- a distribution with the remarkable property that for each $m \in\N$ it can be expressed, locally near $x=0$, as
\begin{equation}
\label{maass9}
\!\! x^{m}P_m\!\left(x\,{\textstyle\frac{d\ }{dx}}\right)\!\psi_{m}(x),\ \text{with $\psi_m$ defined and continuous near the origin};
\end{equation}
here $P_m$ denotes a complex polynomial, whose coefficients depend on $m$ and $\l$. In \cite{inforder} we
introduced the terminology {\it vanishing to infinite order} at $x=0$ for the property (\ref{maass9}) of a distribution
defined on a neighborhood of the origin in $\R$.

To summarize the discussion so far, we have shown that a distribution $\,\tau_\l$ satisfying the automorphy condition
(\ref{maass4}) for $\,\G=SL(2,\Z)$, and additionally the condition $a_0=0$, agrees on $\,\R-\{0\}$ with a distribution that
vanishes to infinite order at $x=0$. Thus either $\,\tau_\l$ itself vanishes to infinite order at $x=0$~-- this is the
meaning of the second condition in (\ref{maass6}), of course -- or else differs from such a distribution by one with
support at the origin. A distribution supported at the origin is a linear combination of the delta function and its derivatives,
and cannot vanish to infinite order at $x=0$ unless it is identically zero. If, contrary to our standing hypothesis, $\,\G$ is a
congruence subgroup of $SL(2,\Z)$, the conditions (\ref{maass6}) must be imposed at each of the cusps of $\,\G$. In that
case the second condition (\ref{maass6}) must also be stated slightly differently.

If $F(x+iy)$ is a Maass form, then so is $F(-x+iy)$. It therefore makes sense to speak of even and odd Maass forms,
i.e., Maass forms such that $F(-x+iy)=\pm F(x+iy)$. Every Maass form can be expressed uniquely as the sum of an
even and an odd Maass form. If $F$ is cuspidal, then so are the even and odd parts. The parity of $F$ affects the
Gamma factors in the functional equation of $L(s,F)$. We shall therefore suppose that $F$, and hence also $\,\tau_\l$,
has a definite parity,
\begin{equation}
\begin{aligned}
\label{maass10}\tau_\l(-x)\ &= \ (-1)^\eta\,\tau_\l(x)\,, \ \ \text{or equivalently}
\\
&a_{-n}\ = \ (-1)^\eta\,a_n\ \ \text{for all $n$}\ \ \ \ (\,\eta\in\Z/2\Z\,)\,.
\end{aligned}
\end{equation}
We also suppose that $F$ is cuspidal, so that $\,\tau_\l$ satisfies
(\ref{maass6}). As one consequence of the parity condition, the
$L$-function
\begin{equation}
\label{maass11}
L(s,F)\ =\  {\sum}_{n\geq 1}\,a_n\,n^{-s+\frac \l 2}\qquad (\,\Re s \gg 0\,)
\end{equation}
completely determines all the $a_n$, and therefore also $\,\tau_\l$
and $F$. We had remarked earlier that $\,\tau_\l$ and $\,\tau_{-\l}$
play essentially symmetric roles unless $\l$ is a negative odd
integer or $\l=0$. Outside of those exceptional cases, the Fourier
coefficients of $\,\tau_\l$ and $\,\tau_{-\l}$ are related by the
factor of propor\-tionality $c_\l |n|^\l$, with $c_\l \neq 0$.
Switching $\,\tau_\l$ and $\,\tau_{-\l}$ has the minor effect of
renormalizing the $L$-function (\ref{maass11}) by the non-zero
constant $c_\l$. It is not difficult to eliminate the remaining ambiguity in
normalizing $L(s,F)$, but we shall not pursue the matter here.

Arguing exactly as in the case of a modular form, we see that the {\it signed Mellin transform}
\begin{equation}
\label{maass11.5} M_\eta(s,\tau_\l) \ = \ \int_\R \tau_\l(x)\, (\sg
x)^\eta\,|x|^{s-1}\,dx
\end{equation}
is a well defined entire holomorphic function. It is legitimate to
substitute the Fourier series (\ref{maass5}) for $\,\tau_\l$ and to
interchange the order of summation and integration, again for the
same reasons as in the case of modular forms, hence
\begin{equation}
\label{maass12}
\begin{aligned}
M_\eta(s,\tau_\l) \ &= \ 2\,{\sum}_{n\geq 1} \,a_n\,n^s\,\int_\R
e(x)\,(\sg x)^\eta\,|x|^{s-1}\,dx
\\
&= \ 2\,G_\eta(s)\,L(s+{\textstyle\frac{\l}2},F)\,,
\end{aligned}
\end{equation}
with
\begin{equation}
\label{maass13}
G_\eta(s)\ = \ \int_\R e(x)\,(\sg x)^\eta\,|x|^{s-1}\,dx \ = \
\begin{cases} \medskip
\,{\textstyle\frac{2\,\G(s)}{(2\pi)^{s}}}\, \cos({\textstyle\frac{\pi s}{2}}) &\text{if}\,\ \eta=0
\\
\,{\textstyle\frac{2\,i\,\G(s)}{(2\pi)^{s}}}\, \sin({\textstyle\frac{\pi s}{2}}) &\text{if}\,\ \eta=1\ ;
\end{cases}
\end{equation}
the explicit formula for $G_\eta(s)$ follows from (\ref{gl2,11}). Since
$M_\eta(s,F)$ is entire, (\rangeref{maass12}{maass13}) show that $L(s,F)$
extends meromorphically to the entire $s$-plane.

The change of variables $x\mapsto -1/x$ in (\ref{maass11.5}), combined with the
transformation rule (\ref{maass7}), gives the functional equation
\begin{equation}
\label{maass14}
M_\eta(s,\tau_\l)\ = \ (-1)^\eta\,M_\eta(1-s-\l,\tau_\l)\,,
\end{equation}
which in turn implies the functional equation
\begin{equation}
\label{maass15}
G_\eta(s-{\textstyle\frac{\l}{2}})\,L(s,F)\ = \ (-1)^\eta\,G_\eta(1-s-{\textstyle\frac{\l}{2}})\,L(1-s,F)
\end{equation}
for $L(s,F)$. Standard Gamma identities establish the equivalence between
Maass' version of the functional equation and (\ref{maass15}).

Though we know that the product $G_\eta(s-\frac \l 2)L(s,F)$ is entire, we
cannot yet conclude that $L(s,F)$ is also entire: unlike $\,\G(s)$, $\,G_\eta(s)\,$
has zeroes. To deal with this problem, we consider the Fourier transform
$\,\widehat\tau_\l\,$ of the tempered distribution $\,\tau_\l$. We use the
normalization $\widehat f(y) = \int_\R f(x)e(-xy)dx$.
Then $e(nx)$, considered as tempered dis\-tribution, has Fourier transform
$\mathcal F e(nx)= \d_{n}(x)=$ Dirac delta function at
$x=n$, and
\begin{equation}
\label{maass16}
\widehat \tau_\l(x)\ = \ {\sum}_{n\neq 0}\,a_n\,\d_{n}(x)\,.
\end{equation}
This distribution visibly vanishes in a neighborhood of the origin, in particular
vanishes to infinite order at $x=0$. According to \cite{inforder}*{theorem 3.19},
the fact that $\,\tau_\l$ vanishes to infinite order at $x=0$ -- cf. (\ref{maass6}) --
implies that $\,\widehat\tau_\l(1/x)$ extends across the origin to a distribution
that vanishes there to infinite order. Since both $\,\widehat\tau_\l(x)$ and
$\,\widehat\tau_\l(1/x)$ have this property, the signed Mellin transform
\begin{equation}
\label{maass17}
\begin{aligned}
M_\eta(s,\widehat \tau_\l)\ &= \ 2\,{\sum}_{n> 0}\,a_n\,n^{s-1}\qquad (\,\Re s \ll 0\,)
\\
&= \ 2\,L(1-s+{\textstyle\frac{\l}{2}},F)
\end{aligned}
\end{equation}
is a well defined, entire holomorphic function. In other words, $L(s,F)$ is entire,
as was to be shown.

The preceding argument essentially applies also to the case of modular forms,
except that one is then dealing with automorphic distributions that are neither
even nor odd, but have only positive Fourier coefficients. In fact, if one considers
modular forms and Maass forms not for $SL(2)$ but for $GL(2)$, a single argument
treats both types of automorphic distributions absolutely uniformly. However, the
case of modular forms is simpler in one respect: the fact that the
$L$-function has no poles requires no special argument.

\section{Pairings of automorphic distributions}\label{pairings}

In the last section we encountered automorphic distributions as distributions on the real line, obtained by a limiting
process. For higher rank groups, it is necessary to take a more abstract point of view, which we shall now explain.

Initially in this section $G$ shall denote a reductive Lie group, $Z^0_G$ the identity component in the center $Z_G$
of $G$, and $\,\G\subset G$ an arithmetically defined subgroup. Note that $G$ acts unitarily on
$L^2(\G\backslash G/Z_G^0)$, via right translation. We consider an irreducible unitary representation
$(\pi,V)$ of $G$ which occurs discretely in $L^2(\G\backslash G/Z_G^0)$,
\begin{equation}
\label{autom1}
j : V \ \hookrightarrow\ L^2(\G\backslash G/Z_G^0)\,.
\end{equation}
Recall the notion of a $C^\infty$ vector for $\pi$\,:\, a vector
$v\in V$ such that $g\mapsto \pi(g)v$ is a $C^\infty$ map from $G$
to the Hilbert space $V$. The space of $C^\infty$ vectors $V^\infty
\subset V$ is dense, $G$-invariant, and gets mapped to
$C^\infty(\G\backslash G/Z_G^0)$ by the embedding (\ref{autom1}).
That makes
\begin{equation}
\label{autom2}
\tau\ = \ \tau_j : V^\infty \ \longrightarrow \ \C\,, \qquad \tau(v)\ = \ jv(e)\,,
\end{equation}
a well defined linear map. It is $\,\G$-invariant because $jv\in C^\infty(\G\backslash G/Z_G^0)$, and is
continuous with respect to the natural topology on $V^\infty$. One should therefore think of $\tau$ as a
$\,\G$-invariant distribution vector for the dual representation $(\pi',V')$ -- i.e., $\tau\in \bigl( (V')^{-\infty}\bigr)^\G$.
Very importantly, $\tau$ determines $j$ completely. Indeed, $j$ is $G$-invariant, so the defining identity
(\ref{autom2}) specifies the value of $jv$, $v\in V^\infty$, not only at the identity, but at any $g\in G$. Since $V^\infty$
is dense in $V$, knowing the effect of $j$ on $V^\infty$ means knowing $j$.

The space $L^2(\G\backslash G/Z_G^0)$ is self-dual, hence if $V$ occurs discretely, so does its dual $V'$. Since
we shall be working primarily with $\tau$, we switch the roles of $V$ and $V'$. From now on,
\begin{equation}
\label{autom3}
\tau \in (V^{-\infty})^\G
\end{equation}
shall denote a $\,\G$-invariant distribution vector corresponding to a discrete embedding
$\,V'\hookrightarrow L^2(\G\backslash G/Z_G^0)$. Not all $\,\G$-invariant distribution vectors correspond to
embeddings into $L^2(\G\backslash G/Z_G^0)$; some correspond to Eisenstein series, and others not even
to those.

The arithmetically defined subgroup $\,\G$ is arithmetic with respect to a particular $\Q$-structure on $G$. If
$P\subset G$ is a parabolic subgroup, defined over $\Q$, with unipotent radical $U$, then $\,\G\cap U$ is a lattice in $U$; in other words,
the quotient $U/(\G\cap U)$ is compact. One calls $\tau\in (V^{-\infty})^\G$ cuspidal if
\begin{equation}
\label{autom4}
\int_{U/(\G\cap U)}\,\pi(u)\tau\,du\ = \ 0 \,,
\end{equation}
for the unipotent radical $U$ of any parabolic subgroup $P$ that is
defined over $\Q$. Since there exist only finitely many
$\,\G$-conjugacy classes of such parabolics, cuspidality amounts to
only finitely many conditions. Essentially by definition, cuspidal
embeddings $\,V'\hookrightarrow L^2(\G\backslash G/Z_G^0)$
correspond to cuspidal distribution vectors $\tau\in
(V^{-\infty})^\G$,\, and conversely every cuspidal automorphic
$\tau$ arises from a cuspidal embedding of $V'$ into
$L^2(\G\backslash G/Z_G^0)$.

To get a handle on $\tau\in (V^{-\infty})^\G$, we realize the space of $C^\infty$ vectors $V^\infty$ as a subspace
$V^\infty\hookrightarrow V_{\l,\d}^\infty$ of the space of $C^\infty$ vectors  for a not-necessarily-unitary principal series representation $(\pi_{\l,\d},V_{\l,\d})$. The Casselman embedding theorem \cite{casselman1} guarantees the
existence of such an embedding. For the moment, we leave the meaning of the subscripts $\l,\d$ undefined. They
are the parameters of the principal series, which we shall explain presently in the relevant cases. A theorem of
Casselman-Wallach \cites{casselman1,wallach} asserts that the inclusion $V^\infty\hookrightarrow V_{\l,\d}^\infty$
extends continuously to an embedding of the space of distribution vectors,
\begin{equation}
\label{autom5}
V^{-\infty}\ \hookrightarrow \ V_{\l,\d}^{-\infty}\,.
\end{equation}
This allows us to consider the automorphic distribution $\tau$ as a distribution vector for a principal series
representation,
\begin{equation}
\label{autom6}
\tau \in \bigl( V_{\l,\d}^{-\infty}\bigr)^\G\,.
\end{equation}
When $G=SL(2,\R)$, cuspidal modular forms correspond to embeddings of discrete series representations into
$L^2(\G\backslash G)$, and cuspidal Maass forms to embeddings of unitary principal series or complementary series representations. The realization of discrete series representations of $SL(2,\R)$ as subrepresentations of principal series representations is very well known, making (\ref{autom6}) quite concrete. For general groups, the Casselman embeddings cannot be described equally explicitly, nor do they need to be unique. Those are not obstacles to using (\ref{autom6}) in
studying $L$-functions. In fact, the non-uniqueness is sometimes helpful in ruling out poles of $L$-functions.

Our tool in studying Rankin-Selberg and related $L$-functions is the pairing of automorphic distributions. In this
paper, we shall only discuss Rankin-Selberg $L$-functions for $GL(2)$ and the exterior square $L$-function for
$GL(4)$. Both involve the pairing of automorphic distributions of $GL(2)$. To minimize notational effort, we shall
work with the group
\begin{equation}
\label{autom7}
\begin{aligned}
G\ = \ PGL(2,\R)\ &\cong \ SL^\pm(2,\R)/\{\pm 1\}
\\
&(\,SL^\pm(2,\R) = \{\,g\in GL(2,\R)\,\mid \, \det g = \pm 1\,\}\,)\,,
\end{aligned}
\end{equation}
rather than $G=GL(2,\R)$, for the remainder of this section. We let $B\subset G$ denote the lower triangular
subgroup. For $\l\in\C$ and $\d\in\Z/2\Z$, we define
\begin{equation}
\label{autom8}
\chi_{\l,\d} : B \rightarrow \ \C^*\,,\ \ \ \chi_{\l,\d}\(\begin{smallmatrix}a&0\\c&d\end{smallmatrix} \) =
(\sg {\textstyle\frac{a}{d}})^\d\,|{\textstyle \f{a}{d}}|^{\frac{\l}{2}} \,.
\end{equation}
The parameterization of the principal series involves a ``$\rho$-shift", i.e., a shift by the half-sum of the positive roots. In
our concrete setting
\begin{equation}
\label{autom8.5}
\rho\ = \ 1 \,,
\end{equation}
and we shall write $\chi_{\l-\rho,\d}$ instead of $\chi_{\l-1,\d}$ to be consistent with the usual notation in the subject.
The space of $C^\infty$ vectors for the principal series representation $\pi_{\l,\d}$ is
\begin{equation}
\label{autom9}
\!\! V_{\l,\d}^\infty   =   \left\{ F \in C^\infty(G)  \mid    F(gb) = \chi_{\l-\rho,\d}(b^{-1} ) F(g) \ \text{for all}\ g\!\in\! G,\, b\!\in\! B\,\right\} ,
\end{equation}
with action
\begin{equation}
\label{autom10}
\bigl(\pi_{\l,\d}(g)F\bigr)(h) \, = \, F(g^{-1}h)    \qquad (\, F\in V_{\l,\d}^\infty\,,\,\ g,h\in G\,) \,.
\end{equation}
Quite analogously
\begin{equation}
\label{autom11}
\! V_{\l,\d}^{-\infty} =   \left\{ \tau\! \in\! C^{-\infty}(G)  \mid    \tau(gb) = \chi_{\l-\rho,\d}(b^{-1} )\, \tau(g) \ \text{for all}\ g\!\in\! G,\, b\!\in \! B\right\}
\end{equation}
is the space of distribution vectors, on which $G$ acts by the same formula as on $V_{\l,\d}^\infty$.

The tautological action of $GL(2,\R)$ on $\R^2$ induces a transitive action of $G=PGL(2,\R)$ on $\RP^1$; in fact
$\RP^1 \cong G/B$, since $B$ is the isotropy subgroup at the line spanned by the second standard basis vector of
$\R^2$. According to the so-called ``fundamental theorem of projective geometry''\!, the action of $G$ on $\RP^1$
induces a simply transitive, faithful action on the set of triples of distinct points in $\RP^1\times\RP^1\times\RP^1$.
Put differently, $G$ has a dense open orbit in
\begin{equation}
\label{autom12}
\RP^1\times\RP^1\times\RP^1 \ \cong \ G/B \times G/B \times G/B\,,
\end{equation}
and can be identified with that dense open orbit once a base point has been chosen. The three matrices
\begin{equation}\label{flagchoice}
    f_1 \ \ = \ \ \(\begin{smallmatrix}  1 & 0 \\ 0 & 1
    \end{smallmatrix}\)\,, \  \  \ f_2 \ \ = \ \ \(\begin{smallmatrix}  1 & 1 \\ 0 & 1
    \end{smallmatrix}\)\,, \ \  \ f_3 \ \ = \ \ \(\begin{smallmatrix}  0 & -1 \\  1 & 0  \end{smallmatrix}\)
\end{equation}
lie in distinct cosets of $B$, so
\begin{equation}
\label{autom14}
G \ \hookrightarrow\ \ G/B \times G/B \times G/B\,,\ \ \ g\mapsto \ \bigl(gf_1B,\,gf_2B,\,gf_3B)\,,
\end{equation}
gives a concrete identification of $G$ with its open orbit in $\RP^1\times\RP^1\times\RP^1$.

Formally at least, the existence of the open orbit can be used to define a $G$-invariant trilinear pairing
\begin{equation}
\label{autom15}
\begin{aligned}
&V_{\l_1,\d_1}^\infty \times V_{\l_2,\d_2}^\infty\times V_{\l_3,\d_3}^\infty \ \longrightarrow \ \C\,,
\\
&\ \ \ (F_1,F_2,F_3) \mapsto  P(F_1,F_2,F_3)\, =_{\text{def}} \int_G
F_1(gf_1)\,F_2(gf_2)\,F_3(gf_3)\,dg \,,
\end{aligned}
\end{equation}
between any three principal series representations. Although the $G$-invariance of the pairing is
obvious from this formula, it is not clear that the integral converges. Before addressing the question of
convergence, we should remark that the ``fundamental theorem of projective geometry'' is field-independent.
The same ideas have been used to construct triple pairings for representations of $PGL(2,\Q_p)$. We
should also point out that a different choice of base points $f_j$ would have the effect of multiplying the
pairing by a non-zero constant.

The question of convergence of the integral (\ref{autom15}) is
most easily understood in terms of the ``unbounded realization" of the principal series, which we discuss
next. The subgroup
\begin{equation}
\label{autom16}
N\ = \ \left\{\,n_x = \ttwo{1}{x}{0}{1}\, \mid \, x\in \R\,\right\} \ \cong \ \R
\end{equation}
of $G$ acts freely on $G/B$, and its image omits only a single point, the coset of
\begin{equation}
\label{autom17}
s\  = \ \ttwo{0}{-1}{1}{0}\,.
\end{equation}
It follows that any $F\in V_{\l,\d}^\infty$ is completely determined by its restriction to $N\cong \R$; the defining
identities (\rangeref{autom8}{autom9}) imply that $\phi_0=$ restriction of $F$ to $\R\,$ is related to $\phi_\infty=$ restriction of
$\pi_{\l,\d}(s)F$ to $\R\,$ by the identity $\phi_\infty(x)= |x|^{\l-1}\phi(-1/x)$. This leads naturally to the identification
\begin{equation}
\label{autom18}
V_{\l,\d}^\infty\ \cong \ \left\{\, \phi\in C^\infty(\R)\,\mid \, |x|^{\l-1}\,\phi(-1/x) \in C^\infty(\R)\,\right\}\,,
\end{equation}
with action
\begin{equation}
\label{autom19}
\begin{aligned}
\bigl(\pi_{\l,\d}(g)\phi\bigr)(x) \, = \, \bigl(\sg (ad-bc)\bigr)^\d\bigl( {\textstyle\frac{|cx+d|}{\sqrt{|ad-bc|}}}\bigr)^{\l-1}\, \phi\bigl({\textstyle\frac{ax+b}{cx+d}}\bigr) \qquad
\\
\text{for}\, \ g^{-1}\ =\ \ttwo{a}{b}{c}{d} \in G\, .
\end{aligned}
\end{equation}
If $(\phi_1,\phi_2,\phi_3)\in \bigl(C^\infty(\R)\bigr)^3$ correspond to $(F_1,F_2,F_3) \in
V_{\l_1,\d_1}^\infty \times V_{\l_2,\d_2}^\infty\times V_{\l_3,\d_3}^\infty$ via the unbounded realization
(\ref{autom18}),
\begin{equation}
\label{autom20}
\begin{aligned}
&\!\!\! P(F_1,F_2,F_3)\, = \, \int_{\R^3} \phi_1(x)\,\phi_2(y)\,\phi_3(z)\,k(x,y,z)\,dx\,dy\,dz\,,\ \ \text{with}
\\
&\ \ \  k(x,y,z)\ = \ \sgn\bigl((x-y)(y-z)(z-x)\bigr)^{\d_1+\d_2+\d_3}\ \times
\\
&\qquad\ \ \  \times \ |x-y|^{\frac{-\l_1-\l_2+\l_3-1}{2}} \, |y-z|^{\frac{\l_1-\l_2-\l_3-1}{2}}\, |x-z|^{\frac{-\l_1+\l_2-\l_3-1}{2}}\,.\!
\end{aligned}
\end{equation}
This can be seen from the explicit form of the isomorphism
(\ref{autom18}), coupled with the definition (\ref{autom9}) of
$V_{\l,\d}^\infty$. We should point out that in the setting of Maass forms, $\d$ plays
the role of the parity $\eta$ in (\ref{maass10}).

Contrary to appearance, the integral (\ref{autom20}) is really an
integral over the compact space $\RP^1\times \RP^1\times\RP^1$: the
integral retains the same general form when one or more of the
coordinates $x,\,y,\,z$ are replaced by their reciprocals; this
follows from the behavior of the $\phi_j$ at $\infty$ specified in
(\ref{autom18}). The convergence of the integral is therefore a
purely local matter. Near points where exactly two of the
coordinates coincide, absolute convergence is guaranteed when the
real part of the corresponding exponent is greater than $-1$. To
analyze the convergence near points of the triple diagonal
$\,\{x=y=z\}$, it helps to ``blow up" the triple diagonal in the
sense of real algebraic geometry -- or equivalently, to use polar
coordinates in the normal directions. One then sees that absolute
convergence requires not only the earlier condition
\begin{equation}
\label{conva}
\Re \bigl( \l_i - \l_j - \l_k \bigr)\ > \ -1 \ \ \ \
\text{if}\ \ i\neq j,\, j\neq k,\, k\neq i\,,
\end{equation}
but also
\begin{equation}
\label{convb}
\Re \bigl( \l_1 + \l_2 + \l_3 \bigr)\ <\ 1\,.
\end{equation}
Both conditions certainly hold when the $V_{\l_i,\d_i}$ belong to
the {\it unitary principal series}, i.e., when all the $\l_j$ are
purely imaginary.

The argument we have sketched establishes the existence of an
invariant trilinear pairing between the spaces of $C^\infty$ vectors
of any three unitary principal series representations. The pairing
is known to be unique up to scaling \cite{oks}. Even when the $\l_i$
are not purely imaginary, one can use (\ref{autom20}) to exhibit an
invariant trilinear pairing by meromorphic continuation. Indeed, for
compactly supported functions of one variable, the functional
$f\mapsto \ \int_\R f(x)|x|^{s-1}dx$ extends meromorphically to
$s\in\C$, with first order poles at the non-positive integers, but
no other poles. As was just argued, the integral kernel in
(\ref{autom20}) can be expressed as $|u|^s$ or $|u|^{s_1}
|v|^{s_2}$, in terms of suitable local coordinates, after blowing up
when necessary. Localizing the problem as before, by means of a
suitable partition of unity, one can therefore assign a meaning to
the integral (\ref{autom20}) for all triples $(\l_1,\l_2,\l_3)\in
\C^3$ outside certain hyperplanes, where the integral has poles.
Even for parameters $(\l_1,\l_2,\l_3)$ on these hyperplanes one can
exhibit an invariant triple pairing by taking residues.

Let us now consider the datum of distribution vectors $\tau_j \in
V_{\l_j,\d_j}^{-\infty}$ for three principal series representations
$V_{\l_j,\d_j}$, $1\leq j \leq 3$. The unbounded realization of the
$V_{\l_j,\d_j}^{-\infty}$ is slightly more complicated than the
$C^\infty$ case (\ref{autom18}): unlike a $C^\infty$ function, a
distribution is not determined by its restriction to a dense open
subset of its domain. The distribution analogue of (\ref{autom18}),
\begin{equation}
\label{autom21}
V_{\l,\d}^{-\infty}\ \cong \ \left\{\,
(\sigma_0,\sigma_\infty)\in \bigl(C^{-\infty}(\R)\bigr)^2\,\mid \,
\sigma_\infty(x)=|x|^{\l-1}\,\sigma_0(-1/x)\,\right\}\,,
\end{equation}
therefore involves a pair of distributions on $\R$ that determine
each other on $\R-\{0\}$. Suppose now that $\tau_j \cong
(\sigma_{j,0},\sigma_{j,\infty})$ via (\ref{autom21}). Then
\begin{equation}
\label{autom22}
\begin{aligned}
\!\!\!\! (x,y,z)&\, \mapsto\, \sigma_{1,0}(x)\,\sigma_{2,0}(y)\,
\sigma_{3,0}(z)\sgn\bigl((x-y)(y-z)(z-x)\bigr)^{\!\d_1+\d_2+\d_3}\,\times
\\
&\times \ |x-y|^{\frac{-\l_1-\l_2+\l_3+1}{2}}\,
|y-z|^{\frac{\l_1-\l_2-\l_3+1}{2}}\,
|x-z|^{\frac{-\l_1+\l_2-\l_3+1}{2}}
\end{aligned}
\end{equation}
extends naturally to a distribution on $\{(x,y,z)\in
\bigl(\RP^1\bigr)^3 \mid x\neq y\neq z\neq x\}$; as one or more of
the coordinates tend to $\infty$, one replaces those coordinates by
the negative of their reciprocals, and simultaneously the
corresponding $\sigma_{j,0}$ by $\sigma_{j,\infty}$. Since
$\{(x,y,z)\in \bigl(\RP^1\bigr)^3 \mid x\neq y\neq z\neq x\}\cong G$
via the identification (\ref{autom14}), we may regard
(\ref{autom22}) as a distribution on $G$. In fact, this distribution
is
\begin{equation}
\label{autom23} \bigl\{\,g\ \mapsto\ \tau_1(gf_1)\, \tau_2(gf_2)\,
\tau_3(gf_3)\,\bigr\}\ \in \ C^{-\infty}(G)\,,
\end{equation}
although the latter description has no immediately obvious meaning
without the steps we have just gone through. The apparent
discrepancy between the signs in the exponents in (\ref{autom20})
and (\ref{autom22}) reflects the fact that
\begin{equation}
\label{autom24}
|x-y|^{-1}\,|y-z|^{-1}\,|z-x|^{-1} \,dx\,dy\,dz\ \cong \ dg\ = \
\text{Haar measure on}\,\ G
\end{equation}
via the identification (\ref{autom14}). Let us formally record the
substance of our discussion:

\begin{observation}
\label{obs1}
For $\tau_j \in V_{\l_j,\d_j}^{-\infty}$\,,\, $1\leq
j\leq 3$\,,
\[
g\ \mapsto\ \tau_1(gf_1)\, \tau_2(gf_2)\, \tau_3(gf_3)
\]
is a well defined distribution on $G$.
\end{observation}

To motivate our result on pairings of automorphic distributions, we
temporarily deviate from our standing assumption that $\,\G\subset
G$ be arithmetically defined; instead we suppose that $\,\G\subset
G$ is a discrete, cocompact subgroup. In that case, if $\tau_j \in
(V_{\l_j,\d_j}^{-\infty})^\G$\!,\, $\,1\leq j\leq 3$\,, are
$\,\G$-invariant distribution vectors, (\ref{autom23}) defines a
distribution on the compact manifold $\,\G\backslash G$. As such, it
can be integrated against the constant function $1$, and
\begin{equation}
\label{autom25}
\int_{\G\backslash G} \tau_1(gf_1)\, \tau_2(gf_2)\,
\tau_3(gf_3)\,dg
\end{equation}
has definite meaning. The value of the integral remains unchanged
when the variable of integration $g$ is replaced by $gh$, for any
particular $h\in G$. Thus, if $\psi\in C^\infty_c(G)$ has total
integral one,
\begin{equation}
\label{autom26}
\begin{aligned}
&\int_{\G\backslash G} \tau_1(gf_1)\, \tau_2(gf_2)\,
\tau_3(gf_3)\,dg\ =
\\
&\ \ \ \ \ = \ \int_G\int_{\G\backslash G} \tau_1(ghf_1)\,
\tau_2(ghf_2)\,\tau_3(ghf_3)\,\psi(h)\, dg\,dh
\\
&\ \ \ \ \ = \ \int_{\G\backslash G} \biggl( \int_G \tau_1(ghf_1)\,
\tau_2(ghf_2)\, \tau_3(ghf_3)\,\psi(h)\,dh\biggr) dg \,.
\end{aligned}
\end{equation}
The implicit use of Fubini's theorem at the second step can be
justified by a partition of unity argument. In short, we have
expressed the integral (\ref{autom25}) as the integral over
$\,\G\backslash G$ of the $\,\G$-invariant function
\begin{equation}
\label{autom27}
g\ \mapsto \ \int_G \tau_1(ghf_1)\, \tau_2(ghf_2)\,
\tau_3(ghf_3)\,\psi(h)\,dh \,.
\end{equation}
This function is smooth, like any convolution of a distribution with
a compactly supported $C^\infty$ function. Note that the integral
(\ref{autom27}) is well defined even for parameters
$(\l_1,\l_2,\l_3)\in \C^3$ which correspond to poles of the integral
(\ref{autom20}).

We now return to our earlier setting, of an arithmetically defined
subgroup $\,\G \subset G = PGL(2,\R)$, specifically a congruence
subgroup
\begin{equation}
\label{autom28}
\G \ \subset\ PGL(2,\Z)\,.
\end{equation}
In this context, the integral (\ref{autom25}) has no obvious
meaning, since we would have to integrate a distribution over the
noncompact manifold $\,\G\backslash G$. The ``smoothed" integral,
however, potentially makes sense: if the integrand
(\ref{autom27}) can be shown to decay rapidly towards the cusps of
$\,\G\backslash G$, it is simply an ordinary, convergent integral.
That is the case, under appropriate hypotheses:

\begin{thm}
\label{thm1} Let $\tau_j \in (V_{\l_j,\d_j}^{-\infty})^\G$\!, $1\leq
j\leq 3$, be $\G$-automorphic distributions, and $\psi \in
C^\infty_c(G)$ a test function, subject to the normalizing condition
\[
\int_G \psi(g)\,dg \ = \ 1\,.
\]
If at least one of the $\tau_j$ is cuspidal, the $\,\G$-invariant
$C^\infty$ function
\[
F(g)\ = \ \int_G \tau_1(ghf_1)\, \tau_2(ghf_2)\,
\tau_3(ghf_3)\,\psi(h)\,dh
\]
decays rapidly along the cusps of $\,\G$; in particular
$\,\int_{\G\backslash G} F(g)\, dg\,$ converges absolutely. This
integral does not depend on the specific choice of $\,\psi$. If, in
addition, one of the $\tau_j$ depends holomorphically on a complex
parameter,
\[
\int_{\G\backslash G} F(g)\,dg\ = \ \int_{\G\backslash G}\int_G
\tau_1(ghf_1)\, \tau_2(ghf_2)\, \tau_3(ghf_3)\,\psi(h)\,dh\,dg
\]
also depends holomorphically on that parameter.
\end{thm}

Why does $F$ decay rapidly? It is not a modular form~-- the Casimir
operator of $G$ does not act on it finitely. Nor does $F$ satisfy
the condition of cuspidality. However, $F$ can be expressed as the
restriction to the diagonal of a modular form in three variables:
\begin{equation}
\label{autom29}
(g_1,g_2,g_3) \ \mapsto\ \int_G \tau_1(g_1hf_1)\,
\tau_2(g_2hf_2)\, \tau_3(g_3hf_3)\,\psi(h)\,dh
\end{equation}
is a $C^\infty$ function on $G\times G\times G$; this follows from
the fact that the cosets $f_jB$ lie in general position. Since
$\tau_j \in (V_{\l_j,\d_j}^{-\infty})^\G$\!, (\ref{autom29}) is a
$\,\G$-invariant eigenfunction of the Casimir operator in each of
the variables separately. It is cuspidal in the variable
corresponding to the cuspidal factor $\tau_j$, hence decays rapidly
in this one direction. It has at worst moderate growth in the other
directions, and therefore decays rapidly when restricted to the
diagonal. The remaining assertions of the lemma are relatively
straightforward.

We shall need a variant of the theorem in the last section, for the
analysis of the exterior square $L$-function for $GL(4)$. Two of the
$\tau_j$ then occur coupled, as a distribution vector for a
principal series representation of $G\times G$, $\,\G$-invariant
only under the diagonal action, not separately. These two $\tau_j$
arise from a single cuspidal automorphic distribution $\tau$ for
$GL(4,\R)$. In this situation the rapid decay of $F$ reflects the
cuspidality of $\tau$.

\section{The Rankin-Selberg $L$-function for $GL(2)$}\label{rsgl2}

The argument we are about to sketch parallels the classical
arguments of Rankin \cite{rankin} and Selberg \cite{selberg}, and of Jacquet \cite{jacquetsequel} in the case of Maass
forms.  We shall pair two automorphic distributions against an Eisenstein series. In our setting, of course,
the Eisenstein series is also an automorphic distribution.

We recall the construction of the distribution Eisenstein series from \cite{extsquare}, specialized to the case of
$G=PGL(2,\R)$. To simplify the discussion, we only work at full level~-- in other words, with
\begin{equation}
\label{eisen1}
\G\ = \ PGL(2,\Z)\ \simeq \ SL^\pm(2,\Z)/\{\,\pm 1\,\}\,.
\end{equation}
We define $\d_\infty \in V_{\nu,0}^{-\infty}$ in terms of the unbounded realization (\ref{autom21}):\,
$\d_\infty$ corresponds to $(\sigma_0,\sigma_\infty)$, with $\sigma_0=0$ and $\sigma_\infty =$ Dirac delta function at $0$.
Then $\pi_{\nu,0}(\g)\d_\infty = \d_\infty$ for all $\g\in\G_\infty = \{\g\in\G\mid \g\infty =\infty\}$. In particular, the series
\begin{equation}
\label{eisen2}
E_\nu \in V_{\nu,0}^{-\infty}\,,\ \ \ E_\nu \ = \ \zeta(\nu+1)\, {\sum}_{\g\in\G/\G_\infty}\, \pi_{\nu,0}(\g)\d_{\infty}\,,
\end{equation}
makes sense at least formally. It is $\G$-invariant by construction. Hence, when we describe $E_\nu$ in terms of the
unbounded realization (\ref{autom21}), it suffices to specify the first member $\sigma_0$ of the pair
$(\sigma_0,\sigma_\infty)$. This allows us to regard $E_\nu$ as a distribution on the real line,
\begin{equation}
\label{eisen3}
E_\nu\ \ \simeq \ \ {\sum}_{p,q\in\Z,\,\,q>0}\,\,\, q^{-\nu-1}\,\d_{p/q}(x)\,.
\end{equation}
To see the equivalence of (\ref{eisen2}) and (\ref{eisen3}), we note that $\d_{p/q}(x)$, with $p,q\in\Z$ relatively prime,
corresponds to the translate of $\d_\infty$ under
\begin{equation}
\label{eisen4}
\ttwo{p}{r}{q}{s}\in \G\,,\ \ \text{with $r,s\in\Z$ chosen so that}\,\ sp-rq=1\,.
\end{equation}
The disappearance of the factor $\zeta(\nu+1)$ in (\ref{eisen3}) reflects the fact that we now sum over all pairs of
integers $p,q$, with $q>0$, not over relatively prime pairs.

The integral of the series (\ref{eisen3}) against a compactly supported test function converges uniformly and
absolutely when $\re \nu > 1$. Hence $E_\nu\in V_{\nu,0}^{-\infty}\,$ is well defined for $\,\re\nu>1\,$, and depends
holomorphically on $\nu$ in this region. The periodic distribution (\ref{eisen3}) has a Fourier expansion,
\begin{equation}
\label{eisen5}
E_\nu \ \simeq \ {\sum}_{n\in\Z}\,\, a_n\,e(nx) \,.
\end{equation}
To calculate the Fourier coefficients, we reinterpret the sum as a distribution on $\R/\Z$.  Then
\begin{equation}
\label{eisen6}
\begin{aligned}
a_n\, &= \, \int_{\R/\Z} e(-nx)\, {\sum}_{p,q\in\Z,\,\,q>0\ \ }\, q^{-\nu-1}\,\d_{p/q}(x)\,dx
\\
&= \, {\sum}_{q>0}\,\,{\sum}_{0\leq p <q}\,\, q^{-\nu-1}\,e(-np/q)\ = \
\begin{cases}
{\sum}_{d|n}\,d^{-\nu} &\text{if}\,\ n\neq 0 \\ \zeta(\nu) & \text{if}\,\ n=0\,.
\end{cases}
\end{aligned}
\end{equation}
The $a_n$, $n\neq 0$, are entire functions of $\nu$, whereas $a_0=\zeta(\nu)$ has a pole at $\nu=1$, so
\begin{equation}
\label{eisen7}
\begin{gathered}
\text{$E_\nu\,$ extends meromorphically to the entire complex plane,}
\\
\text{with a single pole at $\nu=1$, of order one.}
\end{gathered}
\end{equation}
We should remark that $\d_\infty$ is even with respect to the involution $x\mapsto -x$. This is the reason why at
full level there is no Eisenstein series of odd parity~-- i.e., no Eisenstein series in $V_{\nu,1}^{-\infty}$.

The Eisenstein series (\ref{eisen2}) satisfies a functional equation, which relates $E_{-\nu}\in V_{-\nu,0}^{-\infty}$ to
$E_{\nu}\in V_{\nu,0}^{-\infty}$ via the intertwining operator
\begin{equation}
\label{eisen8}
J_\nu :  V_{-\nu,0}^{-\infty}\ \longrightarrow\ V_{\nu,0}^{-\infty}\,.
\end{equation}
On the level of $C^\infty$ vectors, and in terms of the unbounded realization (\ref{autom18}), the operator is given by
the formula
\begin{equation}
\label{eisen9}
\bigl(J_\nu\phi\bigr)(x) \ =  \ \int_\R \phi(y)\,|y-x|^{\nu-1}\,dy   \,.
\end{equation}
Because of the condition on $\phi$ at infinity, this integral has no singularity at $y=\infty$. At $y=x$, the integral
converges when $\re \nu > 0$, but continues meromorphically to the entire complex plane. It is known that
the integral transform (\ref{eisen9}) extends continuously from an operator $\,J_\nu : V_{-\nu,0}^{\infty}
\to V_{\nu,0}^{\infty}$ between the spaces of $C^\infty$ vectors, to the operator (\ref{eisen8}). Alternatively and
equivalently, (\ref{eisen8}) can be defined as the adjoint of $\,J_{\nu} : V_{-\nu,0}^{\infty} \to V_{\nu,0}^{\infty}$,
using the natural duality\footnote{The duality which extends the $G$-invariant pairing $V_{\nu,0}^{\infty} \times
V_{-\nu,0}^{\infty}\to\C$ given by integration over $\R$, in terms of the unbounded realization.} between
$V_{\nu,0}^{\infty}$ and $V_{-\nu,0}^{-\infty}$. Either way one sees that
\begin{equation}
\begin{CD}
\label{eisen10}
V_{-\nu,0}^{-\infty}\,\ni \, e(nx) \  @>{\textstyle{J_\nu}}>> \ G_0(\nu)\,|n|^{-\nu}\,e(nx) \, \in \, V_{\nu,0}^{-\infty}\ \ \ (\,n\neq 0\,)\,.
\end{CD}
\end{equation}
Here $G_0(\nu)$ refers to the Gamma factor described in (\ref{maass13}), and $e(nx)$ is shorthand for the pair
$\bigr(e(nx),|x|^{\mp \nu-1}e(-n/x)\bigl)$~-- cf. (\ref{autom21}); the second member of the pair can be given a definite
meaning even at the origin, using the notion of {\it vanishing to infinite order} that was discussed in section \ref{GL2}.

In view of the relation (\ref{eisen10}), $J_\nu$ maps the Fourier series (\ref{eisen5}) for $E_{-\nu}$ to $G_0(\nu)$ times
the corresponding series for $E_{\nu}$, except possibly for the constant term and a distribution supported at infinity.
However, no non-zero linear combination of a constant function and a distribution supported at infinity can be
$\G$-invariant. This proves
\begin{equation}
\label{eisen11}
J_\nu E_{-\nu}\ =  \ G_0(\nu)\,E_\nu   \,.
\end{equation}
That is the functional equation satisfied by the Eisenstein series. The parameter $\nu$ is natural from the point of view of representation theory. In the eventual application, we shall work with
\begin{equation}
\label{eisen12}
s\ = \ (\nu + 1)/2
\end{equation}
instead. Note that $\nu \mapsto -\nu$ corresponds to $s\mapsto 1-s$.

We now fix two automorphic distributions, either of which may arise from a modular form or a Maass form,
\begin{equation}
\label{eisen13}
\tau_1 \in (V_{\l_1,\d_1}^{-\infty})^\G   \ \ \text{and}\ \ \tau_2 \in (V_{\l_2,\d_2}^{-\infty})^\G\,,
\end{equation}
of which at least one is cuspidal. According to (\ref{eisen7}) and theorem \ref{thm1}, the integral
\begin{equation}
\label{eisen14}
P_\nu^\G(\tau_1,\tau_2,E_\nu) \, =\, \int_{\G\backslash G}\int_G \!
\tau_1(ghf_1)\, \tau_2(ghf_2)\, E_\nu(ghf_3)\, \psi(h)\,dh\,dg
\end{equation}
depends meromorphically on $\nu\in\C$, with a potential first order pole at $\nu=1$ but no other singularities. The
subscript $\nu$ is meant to emphasize the fact that the third argument lies in the space $(V_{\nu,0}^{-\infty})^\G$, and
the superscript $\G$ distinguishes this pairing of $\G$-invariant distribution vectors from the pairing (\ref{autom20})
between spaces of $C^\infty$ vectors.

We shall derive the Rankin-Selberg functional equation from the functional equation (\ref{eisen11}) of the Eisenstein series.
Since the latter involves the intertwining operator, we need to know how $J_\nu$ relates $P_{-\nu}^\G$ to
$P_\nu^\G$. First the analogous statement about the pairing (\ref{autom20}): for $\,F_1 \in V_{\l_1,\d_1}^\infty,\ F_2 \in V_{\l_2,\d_2}^\infty,\  F_3 \in V_{-\nu,0}^\infty\,$,
\begin{equation}
\label{eisen15}
\begin{aligned}
&\!\! P(F_1,F_2,J_\nu F_3) \ =
\\
&=\,(-1)^{\d_1+\d_2} \frac{G_{\d_1+\d_2}\bigl({\textstyle\frac{\l_1 - \l_2 - \nu +1}{2}}\bigr) G_{\d_1+\d_2}  \bigl({\textstyle\frac{-\l_1 + \l_2 - \nu +1}{2}}\bigr)}{G_0(1-\nu)} P(F_1,F_2,F_3)\,.
\end{aligned}
\end{equation}
Note that $P(\dots)$ on the left and the right side of the equality refer to the pairing $V_{\l_1,\d_1}^\infty\times
V_{\l_2,\d_2}^\infty \times V_{\nu,0}^\infty \to \C$, respectively $V_{\l_1,\d_1}^\infty\times
V_{\l_2,\d_2}^\infty \times V_{-\nu,0}^\infty \to \C$\,. The Gamma factors $G_\d(\dots)$ have the same meaning as in
(\ref{maass13}). Since both sides of the equality depend meromorphically on $\nu$, it suffices to establish it for
values of $\nu$ in some non-empty open region. In view of (\ref{autom20}) and (\ref{eisen9}), the assertion
(\ref{eisen15}) reduces to the identity
\begin{equation}
\label{eisen16}
\begin{aligned}
&\!\! \int_\R \frac{\bigl(\sgn (y-t)(t-x) \bigr)^{\d_1+\d_2}}{\bigl(\sgn (y-z)(z-x) \bigr)^{\d_1+\d_2}}\,\, |x-t|^{\a-1} \, |y-t|^{\b-1} \, |z-t|^{-\a-\b}\,dt\ =
\\
&\ = \, (-1)^{\d_1+\d_2}  \frac{G_{\d_1+\d_2}(\a) G_{\d_1+\d_2}\bigl(\b)} {G_0(\a+\b)}\,|x-y|^{\a+\b-1} \,|x-z|^{-\b}\,|y-z|^{-\a} \,,
\end{aligned}
\end{equation}
with $\,\a=(-\l_1+\l_2-\nu+1)/2\,$, $\,\b=(\l_1-\l_2-\nu+1)/2\,$. The integral converges in the region $\,\re \a >0\,,\,\re \b > 0$\,,
$\,\re(\a+\b)<1\,$. The uniqueness of the triple pairing ensures that (\ref{eisen15}) must be correct up to a multiplicative
constant. But then (\ref{eisen16}) must also be correct, except possibly for the specific constant of proportionality. That
constant can be pinned down in a variety of ways; see, for example, \cite{extsquare}*{Lemma 4.32}.

A partition of unity argument shows that the quantities $P^\G_\nu(\tau_1,\tau_2,J_\nu E_{-\nu})$ and
$P^\G_{-\nu} (\tau_1,\tau_2,E_{-\nu})$ are related by the same Gamma factors as the global pairings in (\ref{eisen15}).
Combining this information with (\ref{eisen11}) and the standard Gamma identity $G_\d(\nu)G_\d(1-\nu)=(-1)^\d$, we find
\begin{equation}
\label{eisen17}
\begin{aligned}
&\!\! P^\G_{\nu}(\tau_1,\tau_2,E_{\nu}) \ =
\\
&=\, (-1)^{\d_1+\d_2} G_{\d_1+\d_2}\!\bigl(\!{\textstyle\frac{\l_1 - \l_2 - \nu +1}{2}}\!\bigr)G_{\d_1+\d_2}\!\bigl(\!{\textstyle\frac{-\l_1 + \l_2 - \nu +1}{2}}\!\bigr)  P^\G_{-\nu}(\tau_1,\tau_2,E_{-\nu})\,.\!\!
\end{aligned}
\end{equation}
Once we relate $P^\G_{\nu} (\tau_1,\tau_2,E_{\nu})$ to the Rankin-Selberg $L$-function, this identity will turn out be the
functional equation.

We begin by substituting the expression (\ref{eisen2}) for $E_\nu$ in (\ref{eisen14}). Initially we argue formally; the unfolding
step will be justified later, when we see that the resulting integral converges absolutely:
\begin{equation}
\label{eisen18}
\begin{aligned}
&\!\!\! P_\nu^\G (\tau_1,\tau_2,E_\nu)\, = \, \int_{\G\backslash G}\int_G
\tau_1(ghf_1) \, \tau_2(ghf_2) \, E_\nu(ghf_3)\, \psi(h)\,dh\,dg
\\
&\ =\ \zeta(\nu+1)\! \sum_{\G/\G_\infty} \! \int_{\G\backslash G}  \int_G \!
\tau_1(ghf_1) \, \tau_2(ghf_2)\, \d_\infty(\gamma^{-1}ghf_3)\, \psi(h)\,dh\,dg
\\
&\ =\ \zeta(\nu+1) \int_{\G_\infty\backslash G} \int_G
\tau_1(ghf_1)\, \tau_2(ghf_2)\, \d_\infty(ghf_3)\,\psi(h)\,dh\,dg \,.
\end{aligned}
\end{equation}
The integrand for the outer integral on the right is no longer $\G$-invariant, but it is $(\G\cap N)$-invariant, of
course, and has all the other properties of the integrand in (\ref{eisen14}). Those are the properties used in the proof of
theorem \ref{thm1} to establish rapid decay. In other words, the same argument shows that the integrand in (\ref{eisen18})
decays rapidly in the direction of the cusp. However, $\G_\infty \backslash G$ is not ``compact in the directions
opposite to the cusp", and we still need to argue that the integral converges in those directions as well.

Together with the upper triangular unipotent subgroup $N\subset G$, the two subgroups
\begin{equation}
\label{eisen19}
K\ = \ SO(2)/\{\pm1\}\,,\ \ \ A\ = \ \left\{ \left. \, a_t\, = \, \ttwo{e^t}{0}{0}{e^{-t}} \, \right| \, t\in \R\, \right\}
\end{equation}
determine the Iwasawa decomposition
\begin{equation}
\label{eisen20}
G^0\ = \ NAK
\end{equation}
of the identity component $G^0\simeq SL(2,\R)/\{\pm1\}$ of $G$. Since $\G_\infty$ meets both components of $G$, and since
$\G_\infty \cap G^0 = \G\cap N$, we can make the identification $\G_\infty\backslash G\simeq (\G\cap N)\backslash G^0$.
Hence, and because
\begin{equation}
\label{eisen21}
dg\ = \ e^{-2\rho}\, (a)\, dn\, da\, dk\,,\ \ \ \ \text{with}\ \ \ e^\rho (a_t)\ = \ e^t\,,
\end{equation}
the identity (\ref{eisen18}) can be rewritten as
\begin{equation}
\label{eisen22}
\begin{aligned}
&P_\nu^\G (\tau_1,\tau_2,E_\nu) =
\\
&\ \ \ =\ \zeta(\nu+1)  \int_K \int_A \int_{(\G\cap N)\backslash N} \int_G e^{-2\rho}(a)\,
\tau_1(nakhf_1) \,\tau_2(nakhf_2) \ \times
\\
&\qquad\qquad\qquad\qquad\qquad\qquad\qquad\ \ \ \times\ \d_\infty(nakhf_3)\, \psi(h)\,dh\,dn\,da\,dk \,.
\end{aligned}
\end{equation}
As the $t$ tends to $+\infty$, the point $g=na_tk$ moves towards the cusp. In the opposite direction, as $t\to -\infty$, the
integrand in (\ref{eisen22}) grows at most like a power of $e^{-t}$. To  see this, and to determine the rate of
growth or decay, we temporarily regard the three instances of the argument $\,nak\,$ as independent of each other, as in
the discussion around (\ref{autom29}). In the case of the $\tau_j$, the maximum rate of growth is $e^{(-|\re \l_j| + 1)t}$,
and in the case of $\d_\infty$, it is $e^{(\re\nu +1)t}$, without absolute value sign around $\re\nu$. The reason for the
latter assertion is that we know the behavior of $\d_\infty(g)$ when $g$ is multiplied on the left by any $n\in N$~-- unchanged~-- and when $g$ is multiplied on the left by any $a_t \in A$~-- by the factor $e^{(\re\nu +1)t}$; cf. (\ref{eisen26})
below. In short, the integrand in (\ref{eisen22}) can be made to decay as $t\to -\infty$ by choosing $\re \nu$ large enough.
That makes the integral converge absolutely and justifies the unfolding process.

The smoothing function $\psi\in C_c^\infty(G)$ in theorem \ref{thm1} is arbitrary so far, except for the
normalization $\int_G \psi(g)dg = 1$. We can therefore require $\psi$ to have support in $G^0$, and also  impose
the condition
\begin{equation}
\label{kinv1}
\psi(kg)\ = \ \psi(g)\ \ \text{for all $k\in K$, $g\in G$}\,;
\end{equation}
the latter can be arranged by averaging the original function $\psi$ over $K$. The analogue of (\ref{eisen21}) for the $KAN$
decomposition is $dg = e^{2\rho}(a)\,dk\, da\, dn$. Hence
\begin{equation}
\begin{aligned}
\label{kinv2}
\int_A\int_N e^{2\rho}(a)\,\psi(an)\,dn\,da \ = \ 1\,,\ \ \text{or equivalently}\  \int_A \psi_A(a)\,da \, = \, 1\,,
\\
\text{with}\ \  \    \psi_A(a)\, = \,e^{2\rho}(a)\int_N \psi(an)\,dn\ = \ \int_N \psi(na)\,dn\,,
\end{aligned}
\end{equation}
restates the normalization condition for the $K$-invariant function $\psi$.

We had argued earlier that the function $e(\ell x)$, for $\ell\neq 0$, has a canonical extension~-- now viewed as distribution~--
across infinity. That allows us to regard $e(\ell x)$ as a well defined element of the unbounded model (\ref{autom21}).
We can also make sense of the constant function $1$ as element of the unbounded model for $\,\re \l > 0$, and for other
values of $\l$ by meromorphic continuation. Whether or not $\ell$ equals zero, we let $B_{\ell,\l,\d} \in V_{\l,\d}^{-\infty}$
denote the distribution vector that corresponds to $e(\ell x)$. Then
\begin{equation}
\label{eisen23}
\pi_{\l,\d}(n_x)B_{\ell,\l,\d}\ = \ e(-\ell x)\,B_{\ell,\l,\d}\,,\ \ \ \text{and}\ \ B_{\ell,\l,\d}(n_x)\ = \ e(\ell x)\,.
\end{equation}
The latter equation has meaning since $N\subset G/B$ is open and $B_{\ell,\l,\d}$, like any vector in $V_{\l,\d}^{-\infty}$, transforms according to a character under right translation by elements of $B$. We had assumed that at least one
among $\tau_1$ and $\tau_2$ is cuspidal~-- $\tau_1$, say, for definiteness. Then
\begin{equation}
\label{eisen24}
\tau_1\  = \ {\sum}_{\ell \neq 0}\, a_\ell \,B_{\ell,\l_1,\d_1}\,,\ \ \ \tau_2\  = \ {\sum}_{\ell \in\Z}\, b_\ell \,B_{\ell,\l_2,\d_2} \ + \ \dots
\end{equation}
are the Fourier expansions of $\tau_1$ and $\tau_2$. Here $\,\dots\,$ stands for a vector in $V_{\l_2,\d_2}^{-\infty}$ that
is $N$-invariant and supported on $sB\subset G/B$; recall (\ref{autom17}) for the definition of $s\in G$. The series
for $\tau_1$ has no such singular contribution on $sB$, as was explained in (\ref{maass6}) and the passage that follows it.

In (\ref{eisen22}), the process of averaging over $\G\backslash \G_\infty$ from the left and smoothing from the right
commute. Thus, using the fact that $\d_\infty$ and $\,\dots\,$ in (\ref{eisen24}) are $N$-invariant, we find
\begin{equation}
\label{eisen25}
\begin{aligned}
\!\!\!\!\!\! &\int_{(\G \cap N)\backslash N} \int_G
\tau_1(nakhf_1)\, \tau_2(nakhf_2)\, \d_\infty(nakhf_3) \, \psi(h)\,dh\,dn \ =
\\
&\ =\, \sum_{\ell\neq 0} a_\ell \, b_{-\ell}\! \int_G
B_{\ell,\l_1,\d_1}(akhf_1)\, B_{-\ell,\l_2,\d_2}(akhf_2) \,\d_\infty(akhf_3)\,\psi(h)\,dh
\\
&\ =\, \sum_{\ell\neq 0} a_\ell \, b_{-\ell}\! \int_G
B_{\ell,\l_1,\d_1}(ahf_1)\, B_{-\ell,\l_2,\d_2}(ahf_2) \,\d_\infty(ahf_3)\,\psi(h)\,dh
\\
&\ =\, \sum_{\ell\neq 0} a_\ell \, b_{-\ell}\! \int_G
B_{\ell,\l_1,\d_1}(ah)\, B_{-\ell,\l_2,\d_2}(ahn_1) \,\d_\infty(ahs)\,\psi(h)\,dh\,;
\end{aligned}
\end{equation}
at the second step we have used the $K$-invariance of $\psi$, and at the last step, we have inserted the concrete
values $f_1=e$, $f_2 = n_1$, $f_3 = s$~-- cf. (\ref{flagchoice}) and (\ref{autom16}).

When we substitute (\ref{eisen25}) into (\ref{eisen22}), we can make several simplifications. The expression on the
right in (\ref{eisen25}) no longer depends on the variable $k$, so the integral over $K$ in (\ref{eisen22}) can be
omitted. The distribution $\d_\infty$ is supported on $sB\subset G$. Hence, when the variable $h$ in (\ref{eisen22})
is written as $h=k n \tilde a$, with $k\in K$, $n\in N$, $\tilde a\in A$, and $dh = dk\,dn\,d\tilde a$,
the $k$-integration reduces to evaluation at $k=e$. Since $A$ acts via $e^{2\rho}$ on the cotangent space at $sB \in G/B$,
\begin{equation}
\label{eisen26}
\d_\infty(aks)\,dk\, = \, e^{2\rho}(a)\,\d_\infty(ksa^{-1})\,dk = \, \chi_{\nu+\rho}(a)\,\d_\infty(ks)\,dk \ \ \text{for}\ \ a\in A\,.
\end{equation}
It follows that $\d_\infty(ahs)\,dk=\d_\infty(akn\tilde a s)\,dk=\d_\infty(a ks(s^{-1}ns) \tilde a^{-1})\,dk$ contributes
the factor  $\chi_{\nu+\rho}(a)\,\chi_{\nu-\rho}(\tilde a)=e^{-2\rho}(\tilde a)\,\chi_{\nu+\rho}(a\tilde a)$  when it is integrated
over $K$. Effectively we have replaced the integrals over $h\in G$ in (\ref{eisen22}) and (\ref{eisen25}) by
integrals over $NA$. But the integrand being smoothed in (\ref{eisen25}) is already $N$-invariant. Thus, instead of
smoothing over $G$ with respect to $\psi$, we only need to smooth over $A$ with respect to $\psi_{A}$, as defined in (\ref{kinv2}):
\begin{equation}
\label{eisen27}
\begin{aligned}
&P_\nu^\G (\tau_1,\tau_2,E_\nu)\ =\ \zeta(\nu+1)\,\sum_{\ell\neq 0} \,a_\ell \, b_{-\ell}\! \int_A \int_A \! e^{-2\rho}(a\tilde a)\,\ \times
\\
&\qquad\qquad\qquad \times \ B_{\ell,\l_1,\d_1}(a\tilde a) \,B_{-\ell,\l_2,d_2}(a \tilde a n_1) \, \chi_{\nu+\rho}(a \tilde a)\,\psi_A(\tilde a)\,d\tilde a \,da \,.
\end{aligned}
\end{equation}
We parametrize $a,\tilde a\in A$ as $a=a_t$, $\tilde a = a_{\tilde t}\,$, as in (\ref{eisen19}), with $t, \tilde t\in\R$ and
$da = dt$, $d\tilde a = d\tilde t$. Then, in view of the definition (\ref{autom11}) of $V_{\l,\d}^{-\infty}$ and the
characterization (\ref{eisen23}) of $B_{\ell,\l,\d}$,
\begin{equation}
\label{eisen28}
\begin{aligned}
&B_{\ell,\l_1,\d_1}(a_t a_{\tilde t}) \ = \ e^{(1-\l_1)(t+\tilde t)}\,B_{\ell,\l_1,\d_1}(e)\ = \ e^{(1-\l_1)(t+\tilde t)}\,,
\\
&B_{-\ell,\l_2,\d_2}(a_t \tilde a_{\tilde t} n_1)\ = \ e^{(1-\l_2)(t+\tilde t)}\,B_{-\ell,\l_2,\d_2}(a \tilde a n_1a^{-1}\tilde a^{-1})
\\
&\qquad \qquad\qquad\qquad \! = \ e^{(1-\l_2)(t+\tilde t)}\,e(-\ell \, e^{2(t+\tilde t)})\,,
\\
&\chi_{\nu+\rho}(a_t a_{\tilde t})\ = \ e^{(\nu + 1)(t+\tilde t)}\,,\qquad  e^{-2\rho}(a_t a_{\tilde t})\ = \ e^{-2(t+ \tilde t)}\,.
\end{aligned}
\end{equation}
This leads to the equation
\begin{equation}
\label{eisen29}
\begin{aligned}
&P_\nu^\G (\tau_1,\tau_2,E_\nu)\ =\ \zeta(\nu+1)\,\sum_{\ell\neq 0} \,a_\ell \, b_{-\ell} \ \times
\\
&\qquad\qquad \times \  \int_{\R} \int_\R e^{(\nu+1-\l_1-\l_2)(t+\tilde t)} \, e(-\ell\,e^{2(t+\tilde t)}) \,\psi_A(a_{\tilde t})\,d\tilde t \,dt \,.
\end{aligned}
\end{equation}
To simplify this expression further, we set $x=e^{2t}$, $y=e^{2\tilde t}$, and
\begin{equation}
\label{eisen30}
\psi_A(a_{\tilde t})\ = \ \psi_\R(y) \qquad (\,y=e^{2\tilde t}\,)\,.
\end{equation}
Then $dx = 2 e^{2t} dt$, $dy = 2 e^{2\tilde t} d\tilde t$, and the normalization (\ref{kinv2}) becomes
\begin{equation}
\label{eisen31}
\int_0^\infty \psi_\R(y)\frac{dy}y\ = \ 2\,.
\end{equation}
Putting all the pieces together, we find
\begin{equation}
\label{eisen32}
\begin{aligned}
&\!\!\! P_\nu^\G (\tau_1,\tau_2,E_\nu)\ =
\\
&=\ \frac{\zeta(\nu+1)}4\,\sum_{\ell\neq 0} \,a_\ell \, b_{-\ell} \int_0^\infty\!\! \int_0^\infty \!\! (xy)^{\frac{\nu+1-\l_1-\l_2}{2}} \, e(-\ell\,x\,y) \,\psi_\R(y)\,\frac{dy}y\,{\frac{dx}x} \,.
\end{aligned}
\end{equation}
We know from the derivation of this formula that the integral and the sum must converge for $\re \nu \gg 0$, and indeed
they do. Since the smoothing function $\psi_\R$ has compact support in $(0,\infty)$, the inner integral is the Fourier transform of a compactly supported $C^\infty$ function on $\R$. The resulting function of $x$ is smooth at the origin and decays rapidly at infinity. That makes the outer integral converge, provided $\re \nu$ is large enough. A change of variables then shows that the double integral has order of growth $\,O( |\ell |^{\text{Re} (\l_1 + \l_2 -\nu -1)/2})$, so the sum does converge, again for $\re \nu \gg 0$.

If we regard $e(-\ell x)$, $\ell \neq 0$, not as a function, but as a distribution that vanishes to infinite order at infinity, the
integral $\int_0^\infty e(-\ell x)\, x^{\frac{\nu+1-\l_1-\l_2}{2}}\,dx $ converges for $\re \nu \gg 0$, and the smoothing process in (\ref{eisen32}) becomes unnecessary. Taking this approach, we make the change of variables $x\mapsto x/y$, which splits off the integral (\ref{eisen31}). Hence
\begin{equation}
\label{eisen33}
\begin{aligned}
&P_\nu^\G (\tau_1,\tau_2,E_\nu)\ =
\\
&\ =\ \frac{\zeta(\nu+1)}2\,\sum_{\ell\neq 0} \,a_\ell \, b_{-\ell} \int_0^\infty  x^{\frac{\nu+1-\l_1-\l_2}{2}} \, e(-\ell\,x) \,\frac{dx}x
\\
&\ =\ \frac{\zeta(\nu+1)}2\,\sum_{\ell\neq 0} \,a_\ell \, b_{-\ell}\, |\ell |^{\frac{\l_1+\l_2-\nu-1}{2}} \!\!\int_0^\infty \!\! x^{\frac{\nu+1-\l_1-\l_2}{2}} \, e\bigl(-(\sg \ell)x\bigr) \,\frac{dx}x
\\
&\ =\ \frac{\zeta(2s)}2\,\sum_{\ell\neq 0} \,\,a_\ell \, b_{-\ell}\, |\ell |^{\frac{\l_1+\l_2}{2}-s}  \int_0^\infty  x^{s-\frac{\l_1+\l_2}{2}} \, e\bigl(-(\sg \ell)x\bigr) \,\frac{dx}x\ .
\end{aligned}
\end{equation}
At the last step, we have expressed $\nu$ in terms of $s$, as in (\ref{eisen12}).

By definition, the Rankin-Selberg $L$-function of the pair of automorphic distributions $\tau_1$, $\tau_2$ is
\begin{equation}
\label{eisen34}
L(s,\tau_1\otimes \tau_2)\ = \ \zeta(2s)\,{\sum}_{n>0}\ a_n\,b_n\, n^{\frac{\l_1+\l_2}{2}-s}\,.
\end{equation}
Recall that the Fourier coefficients $a_n$, $b_n$ depend on the choice of the embedding parameter $\l_j$ over $-\l_j\,$.\,
The standard $L$-function (\ref{maass11}), and (\ref{modular2}) in the case of modular forms, with $\l=1-k$, are defined
in terms of the renormalized  coefficients $a_n |n|^{\l/2}$. For the same reason the renormalized coefficients appear in the
Rankin-Selberg $L$-function. To make the connection between (\ref{eisen33}) and the $L$-function, notice that
translation by the matrix
\begin{equation}
\label{eisen35}
r\ = \ \ttwo{-1}{0}{0}{1}
\end{equation}
transforms $\tau_j\in (V_{\l_j,\d_j}^{-\infty})^\G$, realized as $\tau_j(x)$ in terms of the unbounded model, to
$(-1)^{\d_j}\tau_j(-x)$. Since $r\in\G$, that means $\tau_j(-x)=(-1)^{\d_j}\tau_j(x)$, i.e.,
\begin{equation}
\label{eisen36}
a_{-n}\ = \ (-1)^{\d_1}\,a_n\,,\ \ \ b_{-n}\ = \ (-1)^{\d_2}\,b_n \,.
\end{equation}
Hence
\begin{equation}
\label{eisen37}
\begin{aligned}
\zeta(2s)\,{\sum}_{\ell > 0} \,\,a_\ell \, b_{-\ell}\  |\ell |^{\frac{\l_1+\l_2}{2}-s}\ &= \ (-1)^{\d_2}\,L(s,\tau_1\otimes \tau_2)\,,
\\
\zeta(2s)\,{\sum}_{\ell < 0} \,\,a_\ell \, b_{-\ell}\  |\ell |^{\frac{\l_1+\l_2}{2}-s}\ &= \ (-1)^{\d_1}\,L(s,\tau_1\otimes \tau_2)\,.
\end{aligned}
\end{equation}
This allows us to re-write (\ref{eisen33}) as
\begin{equation}
\label{eisen38}
\begin{aligned}
&\!\!2\,P_\nu^\G (\tau_1,\tau_2,E_\nu)\ = \ L(s,\tau_1\otimes \tau_2) \ \times
\\
&\ \times  \left\{ (-1)^{\d_2}   \!\!\! \int_{-\infty}^0 \!\! |x|^{s-\frac{\l_1+\l_2}{2}} \, e(x) \,\frac{dx}x \ +\  (-1)^{\d_1} \!\!\! \int^{\infty}_0 \!\! |x|^{s-\frac{\l_1+\l_2}{2}} \, e(x) \,\frac{dx}x \right\}
\\
&\qquad\qquad\qquad\ \ \ =\  (-1)^{\d_1} \,G_{\d_1+\d_2}(s - {\textstyle\frac{\l_1+\l_2}{2}})\,L(s,\tau_1\otimes \tau_2)\,;
\end{aligned}
\end{equation}
recall (\ref{maass13}), and also the relationship $\nu = 2s-1$ between $\nu$ and $s$.

To complete the proof of the functional equation, we combine (\ref{eisen38}) with (\ref{eisen17}) and appeal to the standard
Gamma identity $G_\d(s)G_\d(1-s)= (-1)^\d$:

\begin{prop}
\label{prop1}
The Rankin-Selberg $L$-function satisfies the functional equation
\[
L(1-s,\tau_1\otimes \tau_2) \ = \,\  {\prod}_{\e_1,\e_2 = \pm 1} \ G_{\d_1+\d_2}(s+\e_1{\textstyle\f{\l_1}{2}}+\e_2{\textstyle\f{\l_2}2})\,L(s,\tau_1\otimes\tau_2)\,.
\]
\end{prop}

We have shown that (\ref{eisen38}) has a holomorphic continuation to $\C-\{1\}$, with at most a simple pole at $s=1$.  Traditionally one states the functional equation and analytic continuation not for the expression in (\ref{eisen38}), but
rather for Langlands' completed $L$-function
\begin{equation}
\label{eisen39}
\L(s,\tau_1\otimes\tau_2) \ = \ L_\infty(s,\tau_1\otimes\tau_2) \, L(s,\tau_1\otimes \tau_2) \,,
\end{equation}
whose ``component at infinity" is a product of Gamma factors that depend on the type of the $\tau_j$. If both $\tau_1$
and $\tau_2$ correspond to Maass forms, then
\begin{equation}
\label{eisen40}
\begin{aligned}
&\!\! \text{Maass case:}\ \ \ L_\infty(s,\tau_1\otimes\tau_2) \, =\,  {\prod}_{\e_1,\e_2=\pm 1}\, \G_\R(s+\e_1{\textstyle\f{\l_1}{2}}+\e_2{\textstyle\f{\l_2}{2}+\eta})\,,
\\
&\qquad \qquad \qquad \qquad \qquad\qquad \,\ \text{with}\, \ \eta\in \{ 0,1\}\,,\ \ \eta \equiv \d_1+\d_2\ \ (\operatorname{mod} 2)\,.
\end{aligned}
\end{equation}
Here $\G_\R$ denotes the Artin $\G$-factor $\pi^{-s/2}\G(s/2)$.  If one of the
$\tau_j$, say $\tau_2$ for definiteness, corresponds to a
holomorphic cusp form of weight $k$, then
\begin{equation}
\label{eisen41}
\!\! \text{mixed case:}\ \ \ L_\infty(s,\tau_1\otimes\tau_2)\, =\, \G_\C(s+{\textstyle\f{\l_1}{2}+\f{k-1}{2}})\,
\G_\C(s-{\textstyle\f{\l_1}{2}}+{\textstyle\f{k-1}{2}})\,,
\end{equation}
where $\G_\C(s)=2(2\pi)^{-s}\G(s)$. Finally, when both $\tau_1$ and $\tau_2$ correspond to holomorphic cusp forms, of weights $k_1$ and $k_2$, respectively,
\begin{equation}
\label{eisen42}
\text{modular forms case:}\ \ \ L_\infty(s,\tau_1\otimes\tau_2)=\G_\C(s+{\textstyle\f{k_1+k_2}{2}}-1) \G_\C(s+{\textstyle\f{|k_1-k_2|}{2}}).
\end{equation}
In all cases, the functional equation of the previous proposition directly implies the equality of $\L(s,\tau_1\otimes \tau_2)$
and $\L(1-s,\tau_1\otimes\tau_2)$, up to a  sign; this follows from standard Gamma identities, in particular the identity
$G_\d(s)G_\d(1-s)=(-1)^{\d}$ and the Legendre duplication formula.

Just as important as the functional equation is the assertion of holomorphy: both $L(s,\tau_1\otimes\tau_2)$ and
$\L(s,\tau_1\otimes\tau_2)$ are holomorphic except for potential first order poles at $s=0$ and $s=1$. For the
uncompleted $L$-function this follows from a classical argument of Jacquet \cite{jacquetsequel}*{Lemma 14.7.5}.
His argument does not require any detailed calculations, and holds in great generality.

Once $L(s,\tau_1\otimes\tau_2)$ is known to be holomorphic on $\C-\{0,1\}$, one can deduce the holomorphy of
$\L(s,\tau_1\otimes\tau_2)$ on $\C-\{0,1\}$ from the results of this section, as follows. Because of the functional equation,
it suffices to rule out poles in the region $\{\re s \geq 1/2, s\neq 1\}$. In effect, we must show that all poles of
$L_\infty (s,\tau_1\otimes\tau_2)$ with $\re s\geq 1/2$ are compensated by zeroes of $L(s,\tau_1\otimes\tau_2)$.
This is an issue only in the Maass case: modular forms have weights at least 2, and the parameter $\l$ of a Maass form
necessarily lies in the region $\{ \,|\re \l | <1/2\,\}$. In the Maass case, only one of the four Gamma factors in (\ref{eisen40})
can have a pole with $\Re s \ge1/2$. Maass forms correspond to irreducible principal series representations, which
involve $\l_j$ and $-\l_j$ symmetrically. We can therefore assume that $\Re{\l_j}\ge 0$, in which case the pole can only
come from the factor $\G_{\R}(s-\f{\l_1+\l_2}{2})$, with $\eta=0$, and must occur at $s=\f{\l_1+\l_2}{2}$. But then $\d_1=\d_2$, and $G_{\d_1+\d_2}(s-\f{\l_1+\l_2}{2})=G_0(s-\f{\l_1+\l_2}{2})$ also has a pole at $s=\f{\l_1+\l_2}{2}$.
We know that (\ref{eisen38}) is holomorphic on $\C-\{0,1\}$, thus forcing $L(s,\tau_1\otimes\tau_2)$ to vanish at
$s=\f{\l_1+\l_2}{2}$, as was to be shown.

\section{Exterior Square on $GL(4)$}
Recall that if $\,F$ is a Hecke eigenform on $GL(n,\Z)\backslash GL(n, \R)$, or more gene\-rally, on the quotient of
$GL(n,\R)$ by a congruence subgroup, the standard $L$-function of $\,F\,$ has an Euler product
\begin{equation}
\label{stnlfngln}
    L(s,F) \ \ = \  \ {\prod}_{p}\ {\prod}_{j\,=\,1}^n
    (1\,-\,\a_{p,j} p^{-s})\i\,.
\end{equation}
The exterior square $L$-function is then defined as an Euler product
\begin{equation}
\label{extsqlfndef}
    L(s,F,Ext^2) \ \ = \  \ {\prod}_p\ L_p(s,F,Ext^2)\,,
\end{equation}
whose factor at any unramified prime $p$ equals
\begin{equation}
\label{extsqlocdef}
    L_p(s,F,Ext^2) \  \ = \ \ {\prod}_{1\le j < k \le n}\ (1\,-\,\a_{p,j}\,\a_{p,k}\,p^{-s})\i\,.
\end{equation}
The appropriate definition of the factors $L_p(s,F,Ext^2)$ corresponding to the finitely many ramified primes is still
a subtle issue. Harris and Taylor recently exhibited local factors for the ramified primes that are consistent with
Langlands functoriality principles, in their proof of the local Langlands conjectures for $GL(n)$.  However, Shahidi
had much earlier given a separate definition, which by all expectations agrees with the one provided by Harris-Taylor,
though the agreement of the two definitions is not obvious. Shahidi furthermore proved that the $L$-function with
his definition of the ramified factors satisfies a functional equation of the type Langlands predicted. Since there can
only be one definition which obeys this functional equation, the potential discrepancy between the Harris-Taylor
and Shahidi definitions poses no problem from the point of view of $L$-functions, though it still is a problem for the
group-theoretic definition of the Langlands conjectures. In any case, an argument which produces the analytic
continuation and functional equation of $L(s,F,Ext^2)$ must give a definition which agrees with Shahidi's.

In our paper \cite{extsquare}, we carry out the {\it archimedean analysis} of the exterior square $L$-function for
$GL(n)$; we establish the holomorphy of the partial $L$-function $L_S(s,F,Ext^2)$ and its completion at infinity $\L_S(s,F,Ext^2)$, in both cases with the factors in (\ref{extsqlfndef}) corresponding to the set $S$ of ramified primes omitted. To keep the discussion simple, we avoid the problem of ramification in the present paper by treating only the full level subgroup
$GL(4,\Z)\subset GL(4,\R)$.

By necessity, the notation in this section will not completely agree with that of the earlier sections; in particular, we now set
\begin{equation}
\label{gl4,1}
\! G\, = \, GL(4,\R)\,,\, \ G_0\, = \, SL^{\pm}(2,\R)\,,\, \ \G\, = \, GL(4,\Z)\,,\, \ \G_0\, = \, SL^{\pm}(2,\Z).
\end{equation}
We shall also work with the subgroups
\begin{equation}
\label{gl4,2}
\begin{aligned}
&G_1\ = \ \left\{ \left. \begin{pmatrix} g_1 & 0 \\ 0 & g_2 \end{pmatrix} \, \right| \, g_1,\, g_2 \in GL(2,\R) \right\}\ \subset \ G\,,
\\
&\qquad \G_1 \ = \ \left\{ \left. \begin{pmatrix} \g & 0 \\ 0 & \g \end{pmatrix} \, \right| \, \g \in GL(2,\Z) \right\}\ \subset \ \G \,,
\\
&\qquad\qquad  U\ = \ \left\{ \left. \begin{pmatrix} 1 & u \\ 0 & 1 \end{pmatrix} \, \right| \, u \in M_{2\times 2}(\R) \right\}\ \subset \ G\,.
\end{aligned}
\end{equation}
Note that $G_1 \simeq GL(2,\R)\times GL(2,\R)$ contains $\G_1 \simeq GL(2,\Z)$, but not as an arithmetic subgroup.

Again we let $B\subset G$ denote the lower triangular Borel subgroup, and we define $B_1 = G_1 \cap B$. Each pair
\begin{equation}
\label{gl4,3}
(\mu,\eta) \in \C^4 \times (\Z/2\Z)^4
\end{equation}
determines a character $\chi_{\mu,\eta} : B \to \C^*$,
\begin{equation}
\label{gl4,4}
\chi_{\mu,\eta} \bigl( a_{i,j} \bigr) \ = \ {\prod}_{1\leq i\leq 4}\ |a_{i,i}|^{\mu_i}\, (\sgn a_{i,i})^{\eta_i}\,,
\end{equation}
and by restriction also a character $\chi_{\mu,\eta} : B_1 \to \C^*$. For $G = GL(4,\R)$,
\begin{equation}
\label{gl4,4.5}
\rho \ = \ \textstyle(\,\frac 32,\,\frac 12,\,-\frac 12,\,-\frac 32\,)
\end{equation}
represents the half sum of the positive roots. In analogy to (\ref{autom11}),
\begin{equation}
\label{gl4,5}
W_{\mu,\eta}^{-\infty} \,  = \,  \left\{ \tau\! \in\! C^{-\infty}(G)  \mid    \tau(gb) = \chi_{\mu-\rho,\eta}(b^{-1} )\, \tau(g) \ \text{for all}\ g \!\in \! G,\, b\!\in\! B\right\}
\end{equation}
is the space of distribution vectors for a generic principal series representation of $G$. Principal series representations of
$G_1 \simeq GL(2,\R)\times GL(2,\R)$ are induced from $B_1$, and hence also parameterized by pairs $(\mu,\eta) \in \C^4 \! \times \! (\Z/2\Z)^4$,
\begin{equation}
\label{gl4,6}
\!\! V_{\mu,\eta}^{-\infty}   =   \left\{ \tau\! \in\! C^{-\infty}(G_1)  \mid    \tau(gb) = \chi_{\mu-\rho,\eta}(b^{-1}\! ) \tau(g) \, \text{for\,all}\, g \!\in\! G_{ 1}, b\!\in\! B_1\!\right\} .\!
\end{equation}
Our current use of the notation $V_{\mu,\eta}^{-\infty}$ is not consistent with (\ref{autom11}). Not only is $G_1$ a product
of two copies of $\,GL(2,\R)$, but the representations we consider need not be trivial on the center of $\,GL(2,\R)$, in contrast
to the situation in section \ref{pairings}, where we considered only automorphic distributions for $PGL(2,\R)$. However,
the $\rho$-shift in (\ref{gl4,6}) {\it is\/} consistent with (\ref{autom11}): the quantity $\rho$ defined in (\ref{gl4,4.5}) restricts to the
corresponding quantities for the two factors of $G_1\simeq {GL(2,\R)\times GL(2,\R)}$.

The arithmetic group $\G$ intersects $U\simeq \R^4$ in a lattice, so $(\G\cap U)\backslash U$ is compact. That makes it
possible to define the operator
\begin{equation}
\label{gl4,7}
\begin{aligned}
&A :  \bigl( W_{\mu,\eta}^{-\infty} \bigr)^\G\  \longrightarrow \   \bigl( V_{\mu,\eta}^{-\infty} \bigr)^{\G_1} \,,
\\
&\qquad A \tau (g)\ = \ \int_{(\G\cap U)\backslash U}  \,\tau(u\,g)\, e(- \tr u)\,du  \qquad (\,g \in G_1\,) \,.
\end{aligned}
\end{equation}
What matters is the fact that the $U\cdot G_1$-orbit of the identity coset in $G/B$ is open. One can therefore
restrict any $\tau \in (W_{\mu,\eta}^{-\infty} \bigr)^\G$ to this open subset, and then further to $G_1$, once the
dependence on the variable $u\in U$ has been smoothed out by taking a single Fourier component. The
restriction to $G_1$ still transforms according to $\chi_{\mu,\eta}^{-1}$ under right translation by elements of
$B_1=G_1 \cap B$. This makes $A\tau$ lie in $V_{\mu,\eta}^{-\infty}$. Conjugation by any $\g\in\G_1$ preserves
the character $u\mapsto e(-\tr u)$ of $U$ and the lattice $\G \cap U$. Since $\G_1 \subset \G$, the $\G$-invariance
of $\tau$ ensures the $\G_1$-invariance of $A\tau$.

We now consider a particular cuspidal $\,\tau\in (W_{\mu,\eta}^{-\infty})^\G$. Since $\G$ contains the center of
$SL^{\pm}(4,\R)$, any such $\tau$ must vanish identically unless
\begin{equation}
\label{center1}
{\sum}_{1\leq j\leq 4}\ \eta_j \, = \,  0\, \ \ \text{in}\, \ \Z/2\Z\,.
\end{equation}
We shall also suppose that
\begin{equation}
\label{center2}
{\sum}_{1\leq j\leq 4}\ \mu_j \, = \,  0 \,.
\end{equation}
This is not a serious restriction: it holds automatically when $\tau$ arises from a discrete summand of
$L^2(\G\backslash G/Z_G^0)$, as in (\ref{autom1}). Even when that is not the case, we can make (\ref{center2}) hold by
twisting $\tau$ with an appropriate character of $Z_G^0$, without destroying the $\G$-invariance.

In section \ref{pairings}, we described the pairing of three $PGL(2,\Z)$-automorphic distributions on $PGL(2,\R)$.
By limiting ourselves to the case of $PGL(2,\R)$ we avoided some notational complications in (\ref{autom20}) and
(\ref{autom21}), without essential loss of generality: in the case of full level, $\,-1 \in GL(2,\Z)$ must act trivially on
any automorphic distribution. In the current setting, we do need the pairing for triples of automorphic distributions on
$GL(2,\R)$. Theorem \ref{thm1} remains correct as stated in this more general situation, provided the integration
is performed over  $SL^{\pm}(2,\Z)\backslash SL^{\pm}(2,\R)$~-- the center of $GL(2,\R)$ is noncompact and
remains noncompact even modulo $GL(2,\Z)$. The statement requires the $\G$-invariance of all three of the
arguments $\tau_j$ of the pairing $P$. Formally, at least, invariance  under the diagonal action of $\G$ on the
three arguments suffices to produce a $\G$-invariant integrand for the outer integral in theorem \ref{thm1}. It is the
proof of rapid decay that forces us to assume $\G$-invariance of each factor. In the present setting, $A\tau$
arises from a cuspidal automorphic distribution $\tau$ on $GL(4,\R)$. It is not difficult to adapt the proof of theorem
\ref{thm1} to this case: after smoothing by some $\psi\in C^\infty_c(G_0)$, the product of $A\tau$ with the Eisenstein
series $E_\nu$ does decay rapidly along the cusp.

We again define the Eisenstein series $E_\nu$ by the formula (\ref{eisen2}), but now summing over $\G_0/(\G_0)_\infty$;
since $\,-1 \in (\G_0)_\infty$, (\rangeref{eisen5}{eisen7}) remain correct. We should remark that the pairing of three
automorphic distributions on $GL(2,\R)$ vanishes identically unless $\,-1\in GL(2,\R)$ acts trivially under the diagonal
action. The parity condition (\ref{center1}) implies that $\,-1$ acts trivially under the diagonal action on $A\tau$. But
$\,-1$ also acts trivially on delta function $\d_\infty$, and hence on the Eisenstein series $E_\nu$. In short, the parity
condition imposed by the action of the center is satisfied in our situation. We have assembled all ingredients to make
sense of
\begin{equation}
\label{gl4,8}
P^{\G_0}_\nu(A\tau,E_\nu) \ = \ \int_{\G_0\backslash G_0}\int_{G_0} \, A\,\tau\ttwo{ghf_1}{0}{0}{ghf_2} E_\nu(ghf_3)\, \psi(h) \, dh\,dg .
\end{equation}
As a function of $\nu$ this is holomorphic, except for a potential first order pole at $\nu=1$. What we said in section
\ref{pairings} about the intertwining operator $J_\nu$ and its interaction with the pairing remains valid, except for the
parity subscripts of the Gamma factors in (\ref{eisen15}) and (\ref{eisen17}), since we now work on $GL(2,\R)$. The
roles of $\l_1$ and $\l_2$ are played by, respectively, $\mu_1-\mu_2$ and $\mu_3-\mu_4$, as can be seen by
comparing the definition (\ref{gl4,6}) of $V^{-\infty}_{\mu,\eta}$ to the definition (\ref{autom11}). Thus, and because of
(\ref{center2}), $\frac{\l_1-\l_2}2$ corresponds to $\mu_1+\mu_4$ and $\frac{\l_2-\l_1}2$ corresponds to $\mu_2+\mu_3$.
This explains the arguments of the Gamma factors in the identity
\begin{equation}
\label{fe1}
\begin{aligned}
&P^{\G_0}_{\nu}(A\tau,E_{\nu}) \, =  \ (-1)^{\eta_2+\eta_3}\ \times
\\
&\ \ \ \times \, G_{\eta_1+\eta_4}\bigl(\mu_1+\mu_4 - {\textstyle\frac{\nu -1}{2}}\bigr)\, G_{\eta_2+\eta_3}
\bigl(\mu_2 + \mu_3 - {\textstyle\frac{\nu -1}{2}}\bigr)\,  P^{\G_0}_{-\nu}(A\tau,E_{-\nu})\,,
\end{aligned}
\end{equation}
which takes the place of (\ref{eisen17}) in the current setting. In the special case when $\eta_1=\eta_2$ and
$\eta_3=\eta_4$ -- i.e, when the action of $G_1\cong GL(2,\R)\times GL(2,\R)$ on $A\tau$ drops to
$PGL(2,\R)\times PGL(2,\R)$ -- (\ref{fe1}) agrees with in (\ref{eisen17}), as it must. In the remaining cases the
identity is deduced from the appropriate variant of (\ref{eisen16}); for details see \cite{extsquare}.

The identity (\ref{fe1}) is the source of the functional equation of the exterior square $L$-function, just as (\ref{eisen17})
was the source of the functional equation for the Rankin-Selberg $L$-function $L(s,\tau_1\otimes \tau_2)$. To make the
connection between the identity (\ref{fe1}) and the exterior square $L$-function, we need to consider the Fourier expansion of
$\tau$ on
\begin{equation}
\label{fe2}
N\ = \ \left\{ \left. n(x,u,v) = \left(\begin{smallmatrix}   1 & x_1 & u_1 & v \\ 0 & 1 & x_2 & u_2 \\ 0 & 0 & 1 & x_3 \\ 0 & 0 & 0 & 1
\end{smallmatrix} \right) \ \right| x \in \R^3,\ u \in \R^2, v \in \R\,\right\} \,.
\end{equation}
Since the $N$-orbit through the identity coset in $G/B$ is open, it is legitimate to restrict $\tau$ to $N$. This restriction is
$(\G\cap N)$-invariant, which allows us to regard $\tau$ as lying in $C^{-\infty}\bigl( (\G\cap N)\backslash N\bigr)$. Every
$(\G\cap N)$-invariant smooth function on $N$, and dually every $(\G\cap N)$-invariant distribution, has a Fourier
expansion with components indexed by~-- roughly speaking~-- the irreducible unitary representations of $N$. For the
one dimensional, or {\it abelian\/}, representations this is literally true, but typically infinite dimensional representation contribute more than once, but finitely often. The {\it non-abelian} Fourier components will turn out not to matter for our purposes. Thus we write
\begin{equation}
\label{fe3}
\tau\bigl(n(x,u,v)\bigr)\ = \ {\sum}_{1\leq j\leq 3}\, a_{n_1,n_2,n_3}\, e(n_j\,x_j) \  + \ \dots \ ,
\end{equation}
with $\,\dots\,$ denoting the sum of the non-abelian Fourier components of $\tau$. The $a_{n_1,n_2,n_3}$ with positive
indices $n_j$ determine all the others:
\begin{equation}
\label{duality0}
a_{\epsilon_1 n_1,\epsilon_2 n_2, \epsilon_3 n_3}\ = \ \epsilon_1^{\eta_1}\,\epsilon_2^{\eta_1+\eta_2}\,\epsilon_3^{\eta_1+\eta_2+\eta_3}\, a_{n_1,n_2,n_3}\ \ \ (\,\epsilon_j \in \{\pm1\}\,)\,.
\end{equation}
Indeed, $\tau$ is invariant under the action of all diagonal matrices with entries $\pm 1$, since $\G$ contains these. Each
of them acts on $N$ by conjugation, which has the effect of reversing the signs of some of the coordinates. One can then
use (\ref{gl4,4}) to determine how the $a_{n_1,n_2,n_3}$ change when the signs of one or more of the indices is flipped.

When $\,\tau$ is a Hecke eigendistribution, the Fourier coefficients $a_{n_1,n_2,n_3}$ are related to the Hecke eigenvalues.
Specifically, $k^{\mu_1+\mu_2} a_{1,k,1}$ is the eigenvalue of the Hecke operator $T_{1,k,1}$. The eigenvalues for Hecke
operators indexed by unramified primes can be expressed in terms of the $\a_{j,p}$ in (\ref{stnlfngln}) \cite{shi}. Jacquet and
Shalika \cite{jacquet}*{\S 2} have used this expression to identify the factors $L_p(s,\tau,Ext^2)$ for unramified primes $p$ in
terms of the Hecke eigenvalues -- in complete gene\-rality for all $n$, not just $n=4$. In the case of $GL(4)$,
\begin{equation}
\label{jl1}
L_p(s,\tau,Ext^2)\ = \ (1-p^{-2s})^{-1} \, {\sum}_{k\geq 0}\,a_{1,p^k,1}\,p^{k(\mu_1+\mu_2-s)} \,.
\end{equation}
At full level, when there are no ramified primes, the Euler product of the local factors for all primes, as in (\ref{extsqlfndef}), expresses the exterior square $L$-function as
\begin{equation}
\label{jl2}
\ L(s,\tau,Ext^2) \ = \ \zeta(2s)\, {\sum}_{n\geq 1}\, a_{1,n,1}\, n^{\mu_1+\mu_2-s}\,.
\end{equation}
One can use this as the definition of the exterior square $L$-function whether or not $\,\tau$ is a Hecke eigendistribution.

\begin{lem} When $s$ and $\nu$ are related by the equation $2s=\nu+1$,
\label{lem1}
\[
P_\nu^{\G_0}(A\tau,E_\nu)\ = \ 2\, (-1)^{\eta_2}\,G_{\eta_1+\eta_2}(s-\mu_1-\mu_2)\,G_{\eta_1+\eta_3}(s-\mu_1-\mu_3)\,L(s,\tau,Ext^2)\,.
\]\end{lem}

Since the proof is lengthy, we shall first deduce the functional equation, which follows from the lemma in combination
with (\ref{fe1}), (\rangeref{center1}{center2}), and the standard Gamma identity $G_\d(s)G_\d(1-s)=(-1)^\d$:\medskip

\begin{prop}
\label{prop2}
$L(1-s,\tau,Ext^2) = \! \underset{1\leq i < j \leq 4}{\prod} \! G_{\eta_i+\eta_j}(s-\mu_i-\mu_j)\,L(s, \tau,Ext^2)$.
\end{prop}

This result is originally due to Kim \cite{kimgl4} and, in the special case when $W_{\mu,\eta}$ belongs to the spherical principal series, to Stade \cite{stade}. We refer the reader to our paper \cite{extsquare} for a discussion of the history
of the exterior square $L$-function for $GL(n)$.

The usual statement of functional equation relates the exterior square $L$-function $L(s,\tau,Ext^2)$ for $GL(n)$ to
that of the dual automorphic distribution $\widetilde\tau$. In our case, with $G=GL(4,\R)$, these two $L$-functions coincide;
that makes it possible to state the functional equation without reference to $\widetilde\tau$.

Just as in the case of the Rankin-Selberg $L$-function for $GL(2)$, Jacquet's general argument implies
that $L(s, \tau,Ext^2)$ is holomorphic, except for possible first order poles at $s=0$ and $s=1$ \cite{extsquare}. The fact
that $P_\nu^{\G_0}(A\tau,E_\nu)$ is holomorphic, together with an analysis of the poles and zeros of the Gamma factors,
establishes the holomorphy of the completed exterior square $L$-function, again with the possible exception of first
order poles at $0$ and $1$. We conclude our paper with the proof of the lemma.
\bigskip

\noindent {\it Proof of Lemma \ref{lem1}.}\,  Recall the notational conventions (\ref{gl4,1}); in particular $G_0^0 = SL(2,\R)$
denotes the identity component of $G_0 = SL^\pm(2,\R)$. We shall suppose that the smoothing function $\psi$ is
supported on $G_0^0$, as we did in  section \ref{rsgl2}. We also impose the $K$-invariance condition (\ref{kinv1})
and define $\psi_A$ as we did in (\ref{kinv2}). In section \ref{pairings} we had pointed out that the expression
(\ref{autom29}) is smooth as function of all three variables. For the same reason
\begin{equation}
\label{gl4,9}
(g_1,g_2,g_3)\ \mapsto \ \int_{G_0^0} A\tau\! \ttwo{g_1 hf_1}{0}{0}{g_2 hf_2}
\d_\infty(g_3 hf_3)\, \psi(h) \, dh
\end{equation}
is a $C^\infty$ function on $G_0\times G_0 \times G_0$. It is also an eigenfunction of the Casimir operator in each
of the three variables, of moderate growth since $\tau$ and $\d_\infty$ are distribution vectors. The cuspidality of
$\tau$ implies that the restriction of this function to the triple diagonal decays rapidly in the cuspidal directions. We
can therefore set $g_1=g_2=g_3=g$ and integrate with respect to $g$ over the quotient $\G_{0,\infty}\!\backslash G_0$,
with $\G_{0,\infty} = \{\g\in \G_0 \mid \g\infty=\infty\}$.

In analogy with (\ref{eisen18}), we insert the definition (\ref{eisen2}) of $E_\nu$ into (\ref{gl4,8}) and unfold: for
$\re \nu \gg 0$,
\begin{equation}
\label{gl4,10}
\!\! P^{\G_0}_\nu(A\tau,E_\nu) \, = \, \zeta(\nu+1) \! \int_{\!\G_{0,\infty}\!\backslash G_0}\!\int_{G_0} \!\!\! A\tau\!\ttwo{ghf_1}{0}{0}{ghf_2}\!
\d_\infty(ghf_3) \psi(h)  dh\,dg\,,
\end{equation}
The justification of this step hinges on two facts. First of all, the
function (\ref{gl4,9}) has moderate growth, as was just pointed out Secondly, we know the behavior of $\d_\infty$
under left translation by elements of $A$. From here on we can justify the unfolding exactly as in section \ref{rsgl2}.
In (\ref{gl4,10}) we can replace $G_0$ in the inner integral by the identity  component $G_0^0$ on which $\psi$ is
supported.  Since $\G_{0,\infty}$ meets both connected components of $G_0$, we can also replace $G_0$ by $G_0^0$
in the outer integral, provided we simultaneously replace $\G_{0,\infty}$ by $\G_{0,\infty}^0 = \G_{0,\infty} \cap G_0^0$.
We parameterize $G_0^0$ by the Iwasawa decomposition $g=n_x ak$~-- recall (\ref{autom16}) and
(\rangeref{eisen20}{eisen21}) . To avoid confusion, we now let $N_0$, $A_0$, $K_0$ denote the subgroups of
$G_0^0= SL(2,\R)$  analogous to $N$, $A$, $K$ in sections \ref{pairings} and \ref{rsgl2}. Note that $\G\cap N_0$
has index $2$ in $\G_{0,\infty}^0$, which also contains $-1$, so $(\G \cap N_0)\backslash N_0 A_0 K_0$ covers
$\G_{0,\infty}^0\backslash G_0^0$ twice. Thus
\begin{equation}
\label{gl4,11}
\begin{aligned}
&\!\! P^{\G_0}_\nu(A\tau,E_\nu) \ = \  2\, \zeta(\nu+1) \ \times
\\
&\, \times \int_{A_0}  \int_0^1 \! \! \int_{G_0^0} e^{-2\rho}(a)\, A\tau\! \ttwo{n_x a hf_1}{0}{0}{n_x a hf_2}
\d_\infty(n_x a hf_3)\, \psi(h) \, dh\,dx\, da\,;
\end{aligned}
\end{equation}
we have legitimately omitted the integration over the Iwasawa component $k$ because $\psi$ is $K$-invariant.

Recall the definition (\ref{gl4,7}) of $A\tau$. It will be convenient to replace $\tau$ by $\tau^0$, defined by the
formula
\begin{equation}
\label{av1}
\tau^0(g)\ = \ \int_{(\G\cap Z_N)\backslash Z_N} \tau(ng)\,dn\ = \ \int_0^1 \tau\left( \left(
\begin{smallmatrix} 1 & 0 & 0 & v \\ 0 & 1 & 0 & 0 \\  0 & 0 & 1 & 0 \\ 0 & 0 & 0 & 1     \end{smallmatrix}
\right) g \right)\,dv\,,
\end{equation}
with $Z_N=$ center of $N$. Then $\tau^0$ is invariant under left translation by elements of $Z_N$, and by elements
of $\G\cap N$. We shall also need to know that
\begin{equation}
\label{av2}
\tau^0(s_{2,3}\,g)\ = \ \tau^0(g)\ \ \ \text{for all $g\in G$, \ \ \ with}\ \ s_{2,3}\ = \ \left(
\begin{smallmatrix} 1 & 0 & 0 & 0 \\ 0 & 0 & 1 & 0 \\  0 & 1 & 0 & 0 \\ 0 & 0 & 0 & 1     \end{smallmatrix}
\right)\,.
\end{equation}
Indeed, $s_{2,3}$ is contained in $\G$ and commutes with the one parameter group over which $\tau$ is averaged to
produce $\tau^0$. The passage from $\tau^0$ to $A\tau$ involves averaging over three more variables,
\begin{equation}
\label{av3}
A\tau(g)\ = \ \int_{\R^3/\Z^3}\, \tau^0\left( \left(\begin{smallmatrix}  1 & 0 & u_1 & 0 \\ 0 & 1 & x_2 & u_2 \\ 0 & 0 & 1 & 0 \\ 0 & 0 & 0 & 1  \end{smallmatrix}\right)\, g\right)\,e(-u_1-u_2)\,dx_2\,du_1\,du_2\,.
\end{equation}
Since
\begin{equation}
\label{av4}
\!\!\!\! \left(\begin{smallmatrix}  1 & 0 & u_1 & 0 \\ 0 & 1 & x_2 & u_2 \\ 0 & 0 & 1 & 0 \\ 0 & 0 & 0 & 1  \end{smallmatrix}\right) \!
\left(\begin{smallmatrix}  1 & x & 0 & 0 \\ 0 & 1 & 0 & 0 \\ 0 & 0 & 1 & x \\ 0 & 0 & 0 & 1  \end{smallmatrix}\right)
\, = \,
\left(\begin{smallmatrix}  1 & x & u_1-xx_2 & u_1x \\ 0 & 1 & 0 & u_2 + xx_2\\ 0 & 0 & 1 & x \\ 0 & 0 & 0 & 1  \end{smallmatrix}\right)\!
\left(\begin{smallmatrix}  1 & 0 & 0 & 0 \\ 0 & 1 & x_2 & 0 \\ 0 & 0 & 1 & 0 \\ 0 & 0 & 0 & 1  \end{smallmatrix}\right) ,
\end{equation}
the equations (\rangeref{gl4,11}{av1}) and (\ref{av3})  imply
\begin{equation}
\label{av5}
\begin{aligned}
&P^{\G_0}_\nu(A\tau,E_\nu) \ = \  2\, \zeta(\nu+1) \   \times
\\
&\ \ \ \ \times\, \int_{\! A_0} \int_0^1 \! \int_{G_0^0}  \int_{\R^3/\Z^3} \!\! \tau^0  \! \left( \left(\begin{smallmatrix}  1 & x & u_1 & 0 \\ 0 & 1 & 0 & u_2 \\ 0 & 0 & 1 & x \\ 0 & 0 & 0 & 1  \end{smallmatrix}\right) \!
\left(\begin{smallmatrix}  1 & 0 & 0 & 0 \\ 0 & 1 & x_2 & 0 \\ 0 & 0 & 1 & 0 \\ 0 & 0 & 0 & 1  \end{smallmatrix}\right)
\!\ttwo{a hf_1}{0}{0}{a hf_2} \right)
\\
&\ \ \ \ \times \  e^{-2\rho}(a)\, e(-u_1-u_2)\, \d_\infty(n_x a hf_3)\, \psi(h)  \,dx_2\,du_1\,du_2\, dh\,dx\, da\,.
\end{aligned}
\end{equation}
We appeal to the invariance of $\tau^0$ under the center of $N$ to justify setting the $(1,4)$-entry of the first matrix in the argument of $\tau^0$ equal to zero.

The variable of integration $x$ occurs three times in (\ref{av5}). Since $\d_\infty$ is $N_0$-invariant, we may as well
drop the factor $n_x$ in its argument. When we omit the integration with respect to $x$ and treat the remaining instances
of $x$ as two separate variables, the integrand~-- after averaging over $\R^3/\Z^3$ and smoothing with respect
to $\psi$\,~-- is a $C^\infty$ function of those two variables; this follows from the fact that (\ref{gl4,9}) is separately smooth
in all three arguments. We can therefore replace the single integral with respect to $x$ by a double integral, provided
we multiply the integrand by the delta function, evaluated on the difference of the two variables. Since
\begin{equation}
\label{av6}
\left(\begin{smallmatrix}  1 & k & 0 & 0 \\ 0 & 1 & 0 & 0 \\ 0 & 0 & 1 & \ell \\ 0 & 0 & 0 & 1  \end{smallmatrix}\right)
\left(\begin{smallmatrix}  1 & x_1 & u_1 & 0 \\ 0 & 1 & 0 & u_2 \\ 0 & 0 & 1 & x_3 \\ 0 & 0 & 0 & 1  \end{smallmatrix}\right)
\ \equiv \
\left(\begin{smallmatrix}  1 & x_1+k & u_1 & 0 \\ 0 & 1 & 0 & u_2 \\ 0 & 0 & 1 & x_3+\ell \\ 0 & 0 & 0 & 1  \end{smallmatrix}\right)
\end{equation}
modulo the center of $N$, the integrand in (\ref{av5}) is separately periodic when the remaining instances of the variable
$x$ are uncoupled. The sum
\begin{equation}
\label{gl4,12}
\d_0(x_1-x_3) \ = \   {\sum}_{\ell\in\Z} \ e\bigl( \ell(x_1-x_3)\bigr)
\end{equation}
represents the ``delta function along the diagonal" in $\R^2/\Z^2$. Thus, in view of what we just said,
\begin{equation}
\label{gl4,13}
\begin{aligned}
&P^{\G_0}_\nu(A\tau,E_\nu) \ = \  2\, \zeta(\nu+1)\ \sum_{\ell\in\Z}\int_{\! A_0} \int_{\R^2/\Z^2} \! \int_{G_0^0}  \int_{\R^3/\Z^3}
e\bigl(\ell(x_1-x_3)\bigr)\ \times
\\
&\qquad\qquad \times \ \tau^0  \! \left( \!\left(\begin{smallmatrix}  1 & x_1 & u_1 & 0 \\ 0 & 1 & 0 & u_2 \\ 0 & 0 & 1 & x_3 \\ 0 & 0 & 0 & 1  \end{smallmatrix}\right) \!
\left(\begin{smallmatrix}  1 & 0 & 0 & 0 \\ 0 & 1 & x_2 & 0 \\ 0 & 0 & 1 & 0 \\ 0 & 0 & 0 & 1  \end{smallmatrix}\right)
\!\ttwo{a hf_1}{0}{0}{a hf_2}\! \right) \, e^{-2\rho}(a)
\\
&\qquad\qquad \times \ e(-u_1-u_2)\, \d_\infty(a hf_3)\, \psi(h)  \,dx_2\,du_1\,du_2\, dh\,dx_1\,dx_3\, da\,.
\end{aligned}
\end{equation}
We use the matrix identity
\begin{equation}
\label{av7}
\!\!\! \left(\begin{smallmatrix}  1 & 0 & 0 & 0 \\ 0 & 1 & \ell & 0 \\ 0 & 0 & 1 & 0 \\ 0 & 0 & 0 & 1  \end{smallmatrix}\right) \!
\left(\begin{smallmatrix}  1 & x_1 & u_1 & 0 \\ 0 & 1 & 0 & u_2\!\!  \\ 0 & 0 & 1 & x_3\!\! \\ 0 & 0 & 0 & 1  \end{smallmatrix}\right) \!
\left(\begin{smallmatrix}  1 & 0 & 0 & 0 \\ 0 & 1 & x_2 & 0 \\ 0 & 0 & 1 & 0 \\ 0 & 0 & 0 & 1  \end{smallmatrix}\right)
=
\left(\begin{smallmatrix}  1 & x_1 & u_1-\ell x_1 & 0 \\ 0 & 1 & 0 &\!\! u_2 + \ell x_3 \! \!\! \\ 0 & 0 & 1 & x_3 \\ 0 & 0 & 0 & 1  \end{smallmatrix}\right) \!
\left(\begin{smallmatrix}  1 & 0 & 0 & 0 \\ 0 & 1 & x_2 + \ell & 0 \\ 0 & 0 & 1 & 0 \\ 0 & 0 & 0 & 1  \end{smallmatrix}\right),
\end{equation}
the $(\G\cap N)$-invariance of $\tau^0$, and the change of variables $u_1 \mapsto u_1+\ell x_1$,
${u_2\mapsto u_2 - \ell x_3}$ to eliminate the factor $e\bigl( \ell(x_1-x_3)\bigr)$ in (\ref{gl4,13}) while simultaneously
replacing $x_2$ by $x_2+\ell$. We then combine the $x_2$-integral over $\{0\leq x_2\leq 1\}$  with the sum over $\ell$
into a single integral over $\R\,$:
\begin{equation}
\label{av8}
\begin{aligned}
&P^{\G_0}_\nu(A\tau,E_\nu) \ = \  2\, \zeta(\nu+1) \   \times
\\
&\ \ \ \ \times\, \int_{\! A_0} \int_{\R^2/\Z^2} \! \int_{G_0^0}  \! \int_\R \int_{\R^2/\Z^2} \tau^0  \! \left(\! \left(\begin{smallmatrix}  1 & x_1 & u_1 & 0 \\ 0 & 1 & 0 & u_2 \\ 0 & 0 & 1 & x_3 \\ 0 & 0 & 0 & 1  \end{smallmatrix}\right) \!
\left(\begin{smallmatrix}  1 & 0 & 0 & 0 \\ 0 & 1 & y & 0 \\ 0 & 0 & 1 & 0 \\ 0 & 0 & 0 & 1  \end{smallmatrix}\right)
\!\ttwo{\! a hf_1 \!}{0}{0}{\! a hf_2\!} \!\right)
\\
&\ \ \ \ \times \  e^{\mathstrut{-2\rho}}(a)\, e(-u_1-u_2)\, \d_\infty(a hf_3)\, \psi(h)  \,du_1\,du_2\, dy \, dh\,dx_1\,dx_3\, da\,.
\end{aligned}
\end{equation}
The symbol $y$ instead of $x_2$ is meant to emphasize the new role of this variable.

Recall the invariance of $\tau^0$ under $s_{2,3}$, as defined in (\ref{av2}). Conjugating $s_{2,3}$ across the first matrix
in the argument of $\tau^0$ has the effect of switching the roles of the $x_i$ and the $u_j$,
\begin{equation}
\label{av9}
\begin{aligned}
&P^{\G_0}_\nu(A\tau,E_\nu) \ = \  2\, \zeta(\nu+1) \   \times
\\
&\ \ \ \ \times\, \int_{\! A_0} \int_{\R^2/\Z^2} \! \int_{G_0^0}  \! \int_\R \int_{\R^2/\Z^2} \tau^0  \! \left(\! \left(\begin{smallmatrix}  1 & x_1 & u_1 & 0 \\ 0 & 1 & 0 & u_2 \\ 0 & 0 & 1 & x_3 \\ 0 & 0 & 0 & 1  \end{smallmatrix}\right) \!
\left(\begin{smallmatrix}  1 & 0 & 0 & 0 \\ 0 & 0 & 1 & 0 \\ 0 & 1 & y & 0 \\ 0 & 0 & 0 & 1  \end{smallmatrix}\right)
\!\ttwo{\! a hf_1 \!}{0}{0}{\! a hf_2\!} \!\right)
\\
&\ \ \ \ \times \  e^{\mathstrut{-2\rho}}(a)\, e(-x_1-x_3)\, \d_\infty(a hf_3)\, \psi(h)  \,dx_1\,dx_3\, dy \, dh\,du_1\,du_2\, da\,.
\end{aligned}
\end{equation}
The congruence
\begin{equation}
\label{av10}
\left(\begin{smallmatrix}  1 & x_1 & u_1 & 0 \\ 0 & 1 & 0 & u_2 \\ 0 & 0 & 1 & x_3 \\ 0 & 0 & 0 & 1  \end{smallmatrix}\right)
\ \equiv \
\left(\begin{smallmatrix}  1 & 0 & u_1 & 0 \\ 0 & 1 & 0 & u_2 \\ 0 & 0 & 1 & 0 \\ 0 & 0 & 0 & 1  \end{smallmatrix}\right) \!
\left(\begin{smallmatrix}  1 & x_1 & 0 & 0 \\ 0 & 1 & 0 & 0 \\ 0 & 0 & 1 & x_3 \\ 0 & 0 & 0 & 1  \end{smallmatrix}\right)
\end{equation}
modulo the center of $N$ implies that we can view the integral with respect to $du_1\,du_2$ as projecting $\tau^0$
to the trivial Fourier components with respect to those two variables, whereas the other integrations operate from the right.
Right translation commutes with projection onto the trivial Fourier components, thus allowing us to shift the integration
with respect to $du_1\,du_2$ all the way to the inside. The passage from $\tau$ to $\tau^0$ already involves a
projection. Together with the $du_1\,du_2$-integral, this gives us the projection
\begin{equation}
\label{av11}
\tau \ \mapsto \ \tau_{\text{abelian}} \,,\ \ \ \tau_{\text{abelian}}(g)\ = \ \int_{(\G \cap [N,N])\backslash N} \tau
\bigl( n\,g\bigr)\,dn\,,
\end{equation}
onto the sum of the abelian Fourier coefficients -- equivalently of invariants for the derived group $[N,N]\subset N$. Thus
(\ref{av10}) reduces to
\begin{equation}
\label{av12}
\begin{aligned}
&P^{\G_0}_\nu(A\tau,E_\nu) \ = \  2\, \zeta(\nu+1) \   \times
\\
&\ \ \ \ \times\, \int_{\! A_0}\! \int_{G_0^0}  \! \int_\R \int_{\R^2/\Z^2} \tau_{\text{abelian}}  \! \left(\! \left(\begin{smallmatrix}  1 & x_1 & 0 & 0 \\ 0 & 1 & 0 & 0 \\ 0 & 0 & 1 & x_3 \\ 0 & 0 & 0 & 1  \end{smallmatrix}\right) \!
\left(\begin{smallmatrix}  1 & 0 & 0 & 0 \\ 0 & 0 & 1 & 0 \\ 0 & 1 & y & 0 \\ 0 & 0 & 0 & 1  \end{smallmatrix}\right)
\!\ttwo{\! a hf_1 \!}{0}{0}{\! a hf_2\!} \!\right)
\\
&\ \ \ \ \times \  e^{-2\rho}(a)\, e(-x_1-x_3)\, \d_\infty(a hf_3)\, \psi(h)  \,dx_1\,dx_3\, dy \, dh\, da\,.
\end{aligned}
\end{equation}

Now we argue as we did in the passage from (\ref{eisen25}) to (\ref{eisen33}). First we substitute $e$, $n_1$, $s$ for
$f_1$, $f_2$, $f_3$ as in (\ref{flagchoice}). We then para\-meterize $h\in G_0^0$ as $h= k\tilde a n_{\tilde x}$, and
observe that the argument
of $\d_\infty$ must lie in $N_0A_0 s \{\pm 1\}$ to give a non-zero contribution. At this point the argument diverges slightly
from our earlier argument, where we worked modulo the center of $SL(2,\R)$. There are three instances of the variable
$h$ in (\ref{av12}). When $h$ is replaced by $(-1)\cdot h$,  $\d_\infty$ remains unchanged, and the other two instances
of $h$ effect a hypothetical sign change of $(-1)^{\eta_1+\eta_2+\eta_3+\eta_4}$ -- hypothetical only since
$\sum_{1\leq j\leq 4} \eta_j = 0$; cf. (\ref{center1}). Thus $k=e$ and $k=-1$ contribute equally, in effect doubling the
factor $2$ in (\ref{av12}). Since
\begin{equation}
\label{av13}
\!\!\! \left(\begin{smallmatrix}  1 & x_1 & 0 & 0 \\ 0 & 1 & 0 & 0 \\ 0 & 0 & 1 & x_3 \\ 0 & 0 & 0 & 1  \end{smallmatrix}\right)  \!
\left(\begin{smallmatrix}  1 & 0 & 0 & 0 \\ 0 & 0 & 1 & 0 \\ 0 & 1 & y & 0 \\ 0 & 0 & 0 & 1  \end{smallmatrix}\right)\!
\left(\begin{smallmatrix}  1 & \tilde x & 0 & 0 \\ 0 & 1 & 0 & 0 \\ 0 & 0 & 1 & \tilde x \\ 0 & 0 & 0 & 1  \end{smallmatrix}\right) \,\equiv \,
\left(\begin{smallmatrix}  1 & -\tilde x y & 0 &  0 \\ 0 & 1 & 0 & 0 \\ 0 & 0 & 1 & \tilde x y\!\! \\ 0 & 0 & 0 & 1  \end{smallmatrix}\right)  \!
\left(\begin{smallmatrix}  1 & x_1 & 0 & 0 \\ 0 & 1 & 0 & 0 \\ 0 & 0 & 1 & x_3 \\ 0 & 0 & 0 & 1  \end{smallmatrix}\right)\!
\left(\begin{smallmatrix}  1 & 0 & 0 & 0 \\ 0 & 0 & 1 & 0 \\ 0 & 1 & y & 0 \\ 0 & 0 & 0 & 1  \end{smallmatrix}\right)
\end{equation}
modulo a left factor lying in $[N,N]$, the variable $n_{\tilde x}$ can simply be absorbed into the $dx_1\,dx_3$-integration.
We can therefore replace $\psi(h)\, dh$ by $\,\psi_A(\tilde a)\, d\tilde a$ and the other instances of $h$ by $\tilde a$, as in (\ref{eisen27}). The smoothing by $\psi$ has now been replaced by smoothing with respect to $\psi_A$,
in the single variable $a$. This reflects the fact that the $A$-direction is the only non-compact direction for the integral
(\ref{av5}), aside from the smoothing integral over $h\in G_0^0$, of course\footnote{The integration with respect to $y\in \R$
in the equivalent integral (\ref{av12}) was obtained by unfolding an integral over $\R/\Z$.}. Just as in section \ref{rsgl2},
the smoothing in the variable $a$ will turn out to be unnecessary when we interpret the integrand -- in effect, a Fourier series in one variable, without constant term -- as a distribution which can be made convergent by integration by parts, under our standing assumption that $\re \nu \gg 0$. To summarize, we can eliminate the integration over $h$ and the factor $\psi(h)$ in (\ref{av12}), provided we double the factor $2$, set $h=e$ in the argument of $\tau^0$, and replace $\d_\infty(a hf_3)$ by  $\chi_{\nu+\rho}(a)$, in ana\-logy to (\ref{eisen26}) and the comment that follows it. Finally we combine the factors $e^{-2\rho}(a)$ and $\chi_{\nu+\rho}(a)$ into the single expression $\chi_{\nu-\rho}(a)$\,:
\begin{equation}
\label{av14}
\begin{aligned}
&P^{\G_0}_\nu(A\tau,E_\nu) \ = \  4\, \zeta(\nu+1) \   \times
\\
&\ \ \ \ \ \ \ \ \times\, \int_{\! A_0}\!  \int_\R \int_{\R^2/\Z^2} \! \tau_{\text{abelian}}  \! \left(\!
\left(\begin{smallmatrix}  1 & x_1 & 0 & 0 \\ 0 & 1 & 0 & 0 \\ 0 & 0 & 1 & x_3 \\ 0 & 0 & 0 & 1  \end{smallmatrix}\right) \!
\left(\begin{smallmatrix}  1 & 0 & 0 & 0 \\ 0 & 0 & 1 & 0 \\ 0 & 1 & y & 0 \\ 0 & 0 & 0 & 1  \end{smallmatrix}\right)
\!\ttwo{{}^{\textstyle a}}{{}^{\textstyle 0}}{{}_{\textstyle 0}}{{}_{\textstyle a n_1}}\! \right)
\\
& \ \ \ \ \ \ \ \ \times \  \chi_{\nu-\rho}(a)\, e(-x_1-x_3)\, dx_1\,dx_3\,dy\, da\,.
\end{aligned}
\end{equation}

For each $n\in (\Z-\{0\})^3$, there exists a unique $B_{n,\mu,\eta} \in W_{\mu,\eta}^{-\infty}$ characterized by the
properties
\begin{equation}
\label{bwgl4,1}
\begin{aligned}
&\pi_{\mu,\eta}\bigl(n(x,u,v)\bigr)B_{n,\mu,\eta}\ = \ e(-n_1x_1 - n_2x_2-n_3x_3)\, B_{n,\mu,\eta}\,,
\\
&\qquad\qquad  B_{n,\mu,\eta}\bigl(n(x,u,v)\bigr)\ = \ e(n_1x_1 + n_2x_2 + n_3x_3)
\end{aligned}
\end{equation}
\cite{chm}; these identities are analogous to (\ref{eisen23}) and use the notation (\ref{fe2}). The $B_{n,\mu,\eta}$
corresponding to different values of $n$ are related by the action of the diagonal subgroup $A\subset G$, but this
need not concern us here. The cuspidality of $\tau$ implies that the Fourier coefficients in (\ref{fe3}) vanish whenever
one or more of the indices are zero. Explicitly,
\begin{equation}
\label{bwgl4,2}
a_n\ \neq \ 0\ \ \Longrightarrow \ \ n\in (\Z-\{0\})^3\,.
\end{equation}
Comparing (\ref{bwgl4,1}) to (\ref{fe2}) and the definition (\ref{av11}) of $\tau_{\text{abelian}}$, one finds
\begin{equation}
\label{bwgl4,3}
\tau_{\text{abelian}} \ = \ {\sum}_{n\in (\Z-\{0\})^3} \ a_n\, B_{n,\mu,\eta} \,.
\end{equation}
The inner integral in (\ref{av14}) picks out the terms in the sum corresponding to $n_1=n_3=1$. Hence
\begin{equation}
\label{bwgl4,4}
\begin{aligned}
&\! P^{\G_0}_\nu(A\tau,E_\nu) \ = \  4\, \zeta(\nu+1) \   \times
\\
&\, \times\ {\sum}_{\ell\neq 0}\,a_{1,\ell,1} \int_{\! A_0}\!  \int_\R  B_{1,\ell,1;\mu,\eta}\! \left(\!
\left(\begin{smallmatrix}  1 & 0 & 0 & 0 \\ 0 & 0 & 1 & 0 \\ 0 & 1 & y & 0 \\ 0 & 0 & 0 & 1  \end{smallmatrix}\right)
\!\ttwo{{}^{\textstyle a}}{{}^{\textstyle 0}}{{}_{\textstyle 0}}{{}_{\textstyle a n_1}}\! \right) \chi_{\nu-\rho}(a)\,dy\, da\,.
\end{aligned}
\end{equation}

We parameterize $A_0$ as in (\ref{eisen19}), $A_0 = \{a_t \mid t\in\R\}$. Then $\chi_{\nu-\rho}(a_t)= e^{(\nu-1)t}$;
cf. (\ref{eisen28}). Conjugating $a_t$ across $n_1$ and using the transformation rule (\ref{gl4,5}) that defines
$W_{\mu,\eta}^{-\infty}$, we can re-write (\ref{bwgl4,4})
as follows:
\begin{equation}
\label{bwgl4,5}
\begin{aligned}
& P^{\G_0}_\nu(A\tau,E_\nu) \ = \  4\, \zeta(\nu+1) \ {\sum}_{\ell\neq 0}\,a_{1,\ell,1} \ \times
\\
&\ \ \ \ \times \ \int_{\R^2}  B_{1,\ell,1;\mu,\eta}\! \left(\!
\left(\begin{smallmatrix}  1 & 0 & 0 & 0 \\ 0 & 0 & 1 & 0 \\ 0 & 1 & y & 0 \\ 0 & 0 & 0 & 1  \end{smallmatrix}\right)
\!\ttwo{{}^{\textstyle 1}}{{}^{\textstyle 0}}{{}_{\textstyle 0}}{{}_{\textstyle a_t n_1a_{-t}}}\! \right) e^{(\nu+1-2\mu_1-2\mu_3)t}\,dy\, dt\,.
\end{aligned}
\end{equation}
The passage from (\ref{bwgl4,4}) to (\ref{bwgl4,5}) also depends on the identity (\ref{center2}), which implies
$(1-\mu_1+\mu_2)+(1-\mu_3+\mu_4) = 2-2(\mu_1+\mu_3)$. Note that $a_t n_1 a_{-t} = n_{e^{2t}}$; cf. (\ref{autom16}).
We appeal to the matrix identity
\begin{equation}
\label{av15}
\left( \begin{smallmatrix} 1 & 0 & 0 & 0 \\ 0 & 0 & 1 &0 \\ 0 & 1 & y & 0 \\ 0 & 0 & 0 & 1 \end{smallmatrix} \right) \!
\left( \begin{smallmatrix} 1 & 0 & 0 & 0 \\ 0 & 1 & 0 & 0 \\ 0 & 0 & 1 & z \\ 0 & 0 & 0 & 1 \end{smallmatrix} \right) \ = \
\left( \begin{smallmatrix} 1 & 0 & 0 & 0 \\ 0 & 1 & 1/y & z \\ 0 & 0 & 1 & yz \\ 0 & 0 & 0 & 1 \end{smallmatrix} \right)\!
\left( \begin{smallmatrix} 1 & 0 & 0 & 0 \\ 0 & -1/y & 0 & 0 \\ 0 & 1 & y & 0 \\ 0 & 0 & 0 & 1 \end{smallmatrix} \right) \,,
\end{equation}
with $z=e^{2t}$, the characterization (\ref{bwgl4,1}) of $B_{n,\mu,\eta}$, and once more to (\ref{gl4,5}), to conclude
\begin{equation}
\label{av16}
\begin{aligned}
&\!\! P^{\G_0}_\nu(A\tau,E_\nu) \ = \  4\,(-1)^{\eta_2} \, \zeta(\nu+1) \ {\sum}_{\ell\neq 0}\,  a_{1,\ell,1} \ \times
\\
&\ \times \, \int_{\R^2}  e\bigl( \ell/y + y\,e^{2t}  \bigr) \,(\sg y)^{\eta_2+\eta_3}\,|y|^{\mu_2-\mu_3-1}\, e^{(\nu+1-2\mu_1-2\mu_3)t}\,dy\, dt\,.
\end{aligned}
\end{equation}
We simplify the integrand by making the change of variables $y \mapsto \ell/y$, followed by the substitution $x=|\ell|\, |y|^{-1}\, e^{2t}$.
Then $dx= 2\,x\, dt$, hence
\begin{equation}
\label{av17}
\begin{aligned}
&\!\! P^{\G_0}_\nu(A\tau,E_\nu) \ = \ 2\, (-1)^{\eta_2} \, \zeta(\nu+1) \ {\sum}_{\ell\neq 0}\,  a_{1,\ell,1} \,|\ell|^{\,\mu_1 + \mu_2-\frac{\nu+1}{2}}\,\ \times
\\
&\ \ \times \ \int_0^\infty\!\!   \int_{\R} \, \frac {e\bigl(y + (\sg \ell y)x \bigr)}{(\sg \ell y)^{\eta_2+\eta_3}} \,\, x^{\frac{\nu - 1}2 -\mu_1-\mu_3}   \,|y|^{\frac{\nu - 1}2 -\mu_1-\mu_2}\,dy\, dx\,.
\end{aligned}
\end{equation}
Recall the definition (\ref{jl2}) of the exterior square $L$-function. We now separate the terms in (\ref{av17})
corresponding to positive and negative values of $\ell$. Appealing to (\ref{center1}) and (\ref{duality0}), we find
\begin{equation}
\label{av19}
\begin{aligned}
&\!\!\! P^{\G_0}_\nu(A\tau,E_\nu) \ = \  2\,L({\textstyle\frac{\nu+1}{2}},\,\tau\,,\,Ext^2)\,\ \times
\\
& \times \, \left\{(-1)^{\eta_2} \! \!\int_0^\infty\! \! \int_{\R}  \frac {e\bigl(y + (\sg y)x \bigr)}{(\sg y)^{\eta_2+\eta_3}} \, |x|^{\frac{\nu - 1}2 -\mu_1-\mu_3}   \,|y|^{\frac{\nu - 1}2 -\mu_1-\mu_2}\,dy\, dx \right.
\\
&\left. + \  (-1)^{\eta_4}\!\! \int_{\!-\infty}^0   \int_{\R} \! \frac {e\bigl(y + (\sg y)x \bigr)}{(\sg y)^{\eta_2+\eta_3}} \, |x|^{\frac{\nu - 1}2 -\mu_1-\mu_3}   \,|y|^{\frac{\nu - 1}2 -\mu_1-\mu_2}\,dy\, dx \right\}.
\end{aligned}
\end{equation}
The factor in curly parentheses equals
\begin{equation}
\label{av19}
\begin{aligned}
&\!\!\! \int_{\R^2}  \frac {(-1)^{\eta_2} \,e(x+y)}{(\sg x)^{\eta_1+\eta_3}\, (\sg y)^{\eta_1+\eta_2}} \, |x|^{\frac{\nu - 1}2 -\mu_1-\mu_3}   \,|y|^{\frac{\nu - 1}2 -\mu_1-\mu_2}\,dy\, dx \ =
\\
&\,= \int_{\R} \! \frac {(-1)^{\mathstrut{\eta_2}} \,e(x)}{(\sg x)^{\eta_1+\eta_3}} \, |x|^{\frac{\nu - 1}2 -\mu_1-\mu_3}   dx \, \times \!
\int_{\R} \! \frac {e(y)}{(\sg y)^{\eta_1+\eta_2}} \,|y|^{\frac{\nu - 1}2 -\mu_1-\mu_2} dy
\\
&\,=\ (-1)^{\mathstrut{\eta_2}}\, G_{\eta_1+\eta_3}({\textstyle\frac{\nu + 1}2} -\mu_1-\mu_3) \, G_{\eta_1+\eta_2}({\textstyle\frac{\nu + 1}2} -\mu_1-\mu_2)\ .
\end{aligned}
\end{equation}
That completes the proof of the lemma. \bx

\vspace{1cm}
\begin{tabular}{lcl}
Stephen D. Miller                    & & Wilfried Schmid \\
Institute of Mathematics           & & Department of Mathematics \\
The Hebrew University          & & Harvard University \\
Jerusalem 91904, Israel                     & & Cambridge, MA 02138 \\
\qquad  \qquad \textsc{and}              & & {\tt schmid@math.harvard.edu}\\
Department of Mathematics \\
Hill Center-Busch Campus  \\
Rutgers, The State University of New Jersey \\
 110 Frelinghuysen Rd \\
 Piscataway, NJ 08854-8019\\
 {\tt miller@math.huji.ac.il}  \\

\end{tabular}

\end{document}